\documentclass[a4paper,12pt,titlepage]{report}

\usepackage[dvips]{graphicx,color}

\usepackage[cp1251]{inputenc}

\usepackage[english,russian]{babel}

\usepackage{amsfonts}
\usepackage{amsmath}
\usepackage{amssymb}
\usepackage{amsthm}
\usepackage{amscd}

\usepackage[dvips]{hyperref}

\usepackage[matrix,arrow,curve]{xy}

\newtheorem{proposition}{Утверждение}
\newtheorem{definition}{Определение}
\newtheorem{assignment}{Упражнения}

\title{Категории. Начальный курс}
\author{Г.В. Кондратьев}

\begin{document}

\maketitle

\begin{abstract}
Теория категорий в России не входит в учебные программы университетов и в тоже время является фундаментальной и необходимой частью математического образования. Учебное пособие рассчитано на первое ознокомление с теорией категорий широкого круга читателей, школьников, студентов, научных работников разных специальностей. 
\end{abstract}

\tableofcontents

\chapter{Категории и функторы}

Категории графически представляют собой ориентированные графы, ребра которых можно перемножать, и это умножение ассоциативно. Но в отличии от теории графов в теории категорий не интересуются такими геометрическими свойствами как, например, число компонентов связности графа или длина минимального пути между вершинами, а переходят на более абстрактную точку зрения, например, существования определенных объектов, построения новых категорий из заданной, взаимоотношения между объектами, категориями, и так далее. Теория категорий предлагает язык достаточно выразительный для многих понятий, конструкций и теорий в математике, который часто делает явными какие-то скрытые стороны теории и позволяет по новому взглянуть на вещи.

\section{Определение категории}
\begin{definition}
{\bf Категория} $\cal K$ состоит из двух непересекающихся классов ${\cal K}_0$ (объекты) и ${\cal K}_1$ (стрелки).
Каждой стрелке $f$ однозначно сопоставляется ее начало $df$ и конец $cf$, которые являются объектами, а каждому объекту $X$ единичная стрелка $1_X$.\\
Стрелки обозначаются $f:X\to Y$ или $X\stackrel {f}{\to }Y$, где $X=df$, $Y=cf$. \\
Если $cf=dg$, как для стрелок $f:X\to Y$ и $g:Y\to Z$, то определена композиция $g\circ f:X\to Z$, являющаяся стрелкой. При этом
\begin{itemize}
\item $(h\circ g)\circ f=h\circ (g\circ f)$ \\
\vskip 0.1cm
(ассоциативность)   
\hskip 1.5cm$\vcenter{\xymatrix{
X \ar[rrr]^-{f} \ar[rrrdd]^(0.3){g\circ f} \ar@<1ex>[dd]^-{(h\circ g)\circ f} \ar@<-1ex>[dd]^-{=}_-{h\circ (g\circ f)} &&& Y \ar[dd]^-{g} \ar[llldd]^(0.3){h\circ g}\\
 &&& \\
W    &&& Z \ar[lll]^-{h}
}}$
\item $1_Y\circ f=f\circ 1_X=f$\\
\vskip 0.1cm
(свойство единицы)
\hskip 1.5cm$\xymatrix{
X \ar[r]^-{1_X} \ar@/_1pc/[rr]_-{f} & X \ar[r]^-{f} \ar@/^1.5pc/[rr]^-{f} & Y \ar[r]^-{1_Y} & Y
}$
\end{itemize} 
\end{definition}
\noindent Как обычно, по умочанию будем предполагать, что если задана даграмма, состоящая из объектов и стрелок, то она {\it коммутативна}, то есть произведения всех стрелок, взятых вдоль разных путей между двумя вершинами, совпадают. Иногда удобно принять соглашение, что каждый объект диаграммы снабжен еще единичной стрелкой. В этом случае произведение стрелок, взятых вдоль пути, начинающемся и заканчивающемся в одной вершине, должно равняться единице, для того чтобы диаграмма была коммутативной.\\
Другое название для стрелок -- {\bf морфизмы}. В дальнейшем будем пользоваться обоими. 

\addvspace{0.2cm}
{\bf Примеры.}
\begin{enumerate}
\item $\mathbf {Set}$, категория множеств и функций между ними.
\item $\mathbf {Set}^*$, категория множеств с отмеченной точкой и функций, сохраняющих отмеченную точку.
\item $\mathbf {Mon}$, категория моноидов и гомоморфизмов. Моноид это множество с одной бинарной ассоциативной операцией и единицей. Гомоморфизм -- функция, сохраняющая операцию и единицу.
\item $\mathbf {Grp}$, категория групп и гомоморфизмов. Группа это моноид, в котором каждый элемент обратим.
\item $\mathbf {Ab}$, категория абелевых групп и гомоморфизмов. Абелева группа это группа с коммутативным сложением.
\item $\mathbf {Vect}$, категория векторных пространств и линейных отображений. Векторное пространство это множество, элементы которого можно складывать и умножать на числа. Линейное отображение это отображение, сохраняющее сложение и умножение на числа. 
\item $\mathbf {Top}$, категория топологических пространств и непрерывных отображений. Топологическое пространство это множество, в котором выделены некоторые открытые множества, представляющие окрестности точек, удовлетворяющие некоторым естественным аксиомам (как, например, что пересечение двух открытых множеств снова открыто, и т.\, д.). Непрерывная функция это функция, у которой прообразы открытых множеств открыты. 
\item $\mathbf {Graph}$, категория (ориентированных) графов и отображений между ними. Граф представляет собой множество вершин, некоторые из которых соединены одним или более ребрами. Отображение графов это отображение, переводящее вершины в вершины, а ребра в ребра с той же ориентацией.
\item $(L,<)$, предупорядоченное множество, рассматриваемое как категория. Именно, скажем, что есть одна стрелка между $x,y\in L$, если $x<y$, и ни одной стрелки в противном случае. В предупорядоченном множестве элементы можно сравнивать друг с другом, при этом всегда выполняется $x<x$ (рефлексивность) и $x<y, \ y<z$ дают $x<z$ (транзитивность).
\item $*\, $, категория, состоящая из одного объекта. В этом случае все стрелки можно перемножать, и это умножение ассоциативно, то есть множество стрелок образует моноид. 
\end{enumerate}

\addvspace{0.2cm}
{\bf Замечание.} Во всех примерах не хватает аксиом в определениях объектов категории. Если читатель не знаком с ними, то он легко сможет их найти в каком-нибудь учебнике алгебры или топологии.

\addvspace{0.2cm}
Введем множество ${\cal K}(X,Y):=\{f\, |\, df=X, cf=Y\}$ и назовем его $Hom$-множеством. $Hom$-множество определено для любых двух объектов категории. 

\begin{definition}
{\bf Подкатегория} ${\cal M}$ категории ${\cal K}$ это категория, у которой ${\cal M}_0\subset {\cal K}_0$, ${\cal M}_1\subset {\cal K}_1$, а закон умножения и единичные стрелки такие же как в ${\cal K}$.
Подкатегория называется {\bf полной}, если для любых объектов $X,Y\in {\cal M}_0\, $ ${\cal M}(X,Y)={\cal K}(X,Y)$.
\end{definition}

Полная подкатегория определена набором своих объектов. Например, если мы возьмем два множества $A$ и $B$ и все функции из $A$ в $A$, из $A$ в $B$, из $B$ в $A$ и из $B$ в $B$, то эти данные задают некоторую полную подкатегорию категории $\mathbf {Set}$.

В заключение параграфа приведем формальное определение двойственной категории, которое упрощает формулировки определений и теорем, где приходится иметь дело с отображениями, меняющими направление стрелок на противоположное.

\begin{definition}
Категория ${\cal K}^{op}$ называется {\bf двойственной} категории ${\cal K}$, если она имеет те же объекты и стрелки, но функции взятия начала и конца стрелки переставлены местами, то есть $d_{{\cal K}^{op}}:=c_{\cal K}$ и $c_{\, {\cal K}^{op}}:=d_{\cal K}$ (началом называется конец, а концом начало стрелки в ${\cal K}$). Соответственно произведение меняется на противоположное $g\circ _{{\cal K}^{op}}f:=f\circ _{\cal K}g$.
\end{definition}

\begin{assignment}
\begin{enumerate}
\item Убедиться, что аксиомы категории (ассоциативность и единица) выполняются для всех примеров.
\item Доказать, что $\mathbf {Ab}\subset \mathbf {Grp}$ и $\mathbf {Grp}\subset \mathbf {Mon}$ полные подкатегории.
\item Показать, что каждому моноиду соответсвует категория с одним объектом.
\item Убедиться, что никакого содержательного смысла стрелкам в $\mathbf {Set}$ с противоположными началом и концом не может быть поставлено.
\end{enumerate}
\end{assignment}

\section{Специальные объекты и стрелки}

\begin{definition}
Стрелка $f:X\to Y$ называется 
\begin{itemize}
\item {\bf монострелкой}, если на нее можно сокращать слева: $f\circ g=f\circ h$  $\Rightarrow $ $g=h$, 
\item {\bf эпистрелкой}, если на нее можно сокращать справа: $g\circ f=h\circ f$ $\Rightarrow $ $g=h$,
\item {\bf изострелкой}, если существует обратная стрелка: $\exists f^{-1}:Y\to X$ такая что $f^{-1}\circ f=1_X$, $f\circ f^{-1}=1_Y$.
\end{itemize}
\end{definition}

Если существует (хотя бы одна) изострелка между объектами $X$ и $Y$, то они называются изоморфными, и это обозначается $X\simeq Y$. 

\begin{assignment}
\begin{enumerate}
\item Показать, что в $\mathbf {Set}$ монострелки соответсвуют инъективным отображениям (инъективные отображения разделяют элементы: $x\ne y$ $\Rightarrow $ $f(x)\ne f(y)$), эпистрелки сюръективным отображениям (отображениям на), а изострелки биекциям.
\item Доказать, что изострелка всегда моно- и эпи-. Показать, что обратное неверно на примере категории предпорядка $(L,<)$.
\item Привести пример изострелок в категории.
\item Доказать, что произведение монострелок есть монострелка, произведение эпистрелок есть эпистрелка, произведение изострелок есть изострелка.
\item Доказать, что отношение изоморфизма является отношением эквивалентности (то есть рефлексивно $X\simeq X$, симметрично $X\simeq Y$ $\Rightarrow $ $Y\simeq X$ и транзитивно $X\simeq Y$ и $Y\simeq Z$ $\Rightarrow $ $X\simeq Z$).
\end{enumerate}
\end{assignment}

\begin{definition}
Объект $\mathbf {0}$ называется {\bf начальным}, если существует ровно одна стрелка из него в любой другой объект категории: $\forall X\in {\cal K}_0$ $\exists !f:\mathbf {0}\to X$. Если, наоборот, существует только одна стрелка из любого объекта в данный, то он называется {\bf конечным} и обозначается $\mathbf {1}$. Объект, являющийся одновременно начальным и конечным называется {\bf нулевым} и обозначается также нулем $0$.
\end{definition}

В категории может не существовать ни одного из этих объектов, или существовать, допустим, только один из них.

\begin{proposition}
Начальный объект $\mathbf {0}$ определен однозначно с точностью до изоморфизма, то есть если имеется другой начальный объект $\mathbf {0}'$, тогда существует (единственный) изоморфизм между ними $\mathbf {0}\stackrel {\exists \, !}{\simeq }\mathbf {0}'$. То же самое справедливо для конечного и нулевого объектов.
\end{proposition}

\begin{assignment}
\begin{enumerate}
\item Доказать приведенное утверждение.
\item Найти начальный и конечный объекты в категориях множеств $\mathbf {Set}$ и множеств с отмеченной точкой $\mathbf {Set}^*$.
\item Найти нулевые объекты в категориях групп $\mathbf {Grp}$ и векторных пространств $\mathbf {Vect}$.
\end{enumerate}
\end{assignment}

\section{Функторы}
Функторы играют роль отображений между категориями, сохраняющих категорную структуру.
\begin{definition}
{\bf Функтор} $F:{\cal K}\to {\cal L}$ это отображение, переводящее объекты в объекты, а стрелки в стрелки, то есть
если $X\in {\cal K}_0$, то $FX\in {\cal L}_0$, и если $f\in {\cal K}_1$, то $Ff\in {\cal L}_1$. При этом
\begin{itemize}
\item $F(df)=d(Ff)$, $F(cf)=c(Ff)$ (сохранение начала и конца стрелки),
\item $F(g\circ f)=Fg\circ Ff$ (сохранение умножения для стрелок, для которых оно определено),
\item $F(1_X)=1_{FX}$ (сохранение единичных стрелок).
\end{itemize}
\end{definition}

Другими словами $F(f:X\to Y)=Ff:FX\to FY$, $F(X\stackrel {f}{\to }Y\stackrel {g}{\to }Z)=FX\stackrel {Ff}{\to }FY\stackrel {Fg}{\to }FZ$ и $F(X\stackrel {1_X}{\to }X)=FX\stackrel {1_{FX}}{\to }FX$.

\addvspace{0.2cm}
{\bf Примеры.}
\begin{enumerate}
\item $1_{\cal K}:{\cal K}\to {\cal K}$, {\bf единичный функтор}, который назначает каждому объекту категории сам этот объект и каждой стрелке саму эту стрелку.
\item $U:\mathbf {Grp}\to \mathbf {Set}:
\left\{
\begin{array}{cr}
G\mapsto |G| & (\text{на объектах})\\
(G\stackrel {f}{\to }G')\mapsto (|G|\stackrel {|f|}{\to }|G'|) & (\text{на стрелках})
\end{array}
\right.
$\\
\vskip 0.05cm
где символ, окруженный вертикальными линиями, означает множество или функцию, рассматриваемую без какой-либо дополнительной структуры. Функтор $U$ называется {\bf забывающим}. Он забывает групповую операцию на множествах и то, что функции между группами, гомоморфизмы.\\
Забывающие функторы связываются обычно с каждой категорией структурированных множеств. Так существуют $\mathbf {Top}\to \mathbf {Set}$, $\mathbf {Vect}\to \mathbf {Ab}$, $\mathbf {Vect} \to \mathbf {Set}$, и так далее.
\item {\bf Постоянный функтор} $\Delta _C:{\cal K}\to {\cal L}$ принимает на всех объектах значение $C\in {\cal L}_0$ и на всех морфизмах значение $1_C\in {\cal L}_1$.
\item Неубывающая функция $f:(L,<)\to (L',<)$ из одного предупорядоченного множества в другое является функтором.
\item $Hom$-функтор, определенный объектом $K\in {\cal K}_0$, задается как\\ 
$Hom(K,-):{\cal K}\to \mathbf {Set}:$\\
\vskip 0.05cm
$\left\{
\begin{array}{cr}
\hskip -1.5cmX\mapsto Hom(K,X) & (\text{на объектах})\\
(X\stackrel {f}{\to }Y)\mapsto (Hom(K,X)\stackrel {f^*}{\to }Hom(K,Y))& (\text{на стрелках})
\end{array}
\right.
$ \\
\vskip 0.05cm
где $f^*(g):=f\circ g$.
\end{enumerate} 

Функтор $F:{\cal K}\to {\cal L}$ индуцирует для каждой пары объектов $X,Y\in {\cal K}_0$ отображение $F_{X,Y}:{\cal K}(X,Y)\to {\cal L}(FX,FY)$. Если для каждой пары объектов это отображение инъективно (переводит разные элементы в разные), то функтор называется {\bf инъективным}, если сюръективно (отображение на), то {\bf сюръективным} или {\bf полным}.

\begin{proposition}
Каждый функтор сохраняет изоморфизмы, то есть если $f$ есть изоморфизм, то $Ff$ тоже изоморфизм.
\end{proposition}

Это простое предложение часто является ключевым для доказательства того, что два объекта неизоморфны (именно, если удается построить функтор в категорию, где образы этих объектов неизоморфны). Например, группа состоящая из двух эментов не может быть изоморфна группе, состоящей из трех элементов, так как (применяем забывающий функтор) множества несущие групповую операцию этих групп не находятся в биективном соответствии.  

\begin{assignment}
\begin{enumerate}
\item Доказать приведенное утверждение.
\item Показать, что полный функтор может не сохранять моно- и эпи- стрелки.
\item Доказать, что неубывающая функция между предупорядоченными множествами является функтором.
\item Доказать, что вложение подкатегории (в категорию) является функтором. Этот функтор полный, если и только если подкатегория полная.
\item Убедиться, что $\mathbf {Set}$ является объектом, а также подкатегорией категории $\mathbf {Cat}$. 
\item Доказать, что $Hom$-функтор является функтором.
\end{enumerate}
\end{assignment}

\begin{proposition}
Существует категория $\mathbf {Cat}$, объекты которой категории, а стрелки -- функторы между ними.
\end{proposition}

\noindent Это предложение позволяет применять теорию категорий к самой себе.

\addvspace{0.2cm}
Часто встречаются отображения между категориями, которые удовлетворяют всем свойствам функтора, но меняют местами начало и конец стрелок. Такие отображения называют {\bf контравариантными} функторами или {\bf предпучками}. Обычные функторы называются тогда {\bf ковариантными}.

\addvspace{0.2cm}
{\bf Примеры.}
\begin{enumerate}
\item {\bf Тождественный контравариантный функтор} $I:{\cal K}\to {\cal K}^{op}$ действует тождественно на объектах и стрелках, но начало и конец стрелок в ${\cal K}^{op}$ задаются наоборот $c_{op}:=d$, $d_{op}:=c$.\\
Функтор $I$ имеет обратный $I^{-1}:{\cal K}^{op}\to {\cal K}$, который действует в точности по такому же правилу как и $I$. Поэтому, любому контравариантному функтору $F:{\cal K}\to {\cal L}$ биективно соответсвует ковариантный функтор $F\circ I^{-1}:{\cal K}^{op}\to {\cal L}$, который обозначают той же буквой $F$, и который определен на двойственной категории ${\cal K}^{op}$ (другими словами, контравариантный функтор $F:{\cal K}\to {\cal L}$ и ковариантный функтор $F:{\cal K}^{op}\to {\cal L}$ считаются одинаковыми).  
\item {\bf Контравариантный} $Hom$-функтор, определенный объектом $K\in {\cal K}_0$, задается как
$Hom(-,K):{\cal K}\to \mathbf {Set}:$\\
\vskip 0.05cm
$\left\{
\begin{array}{cr}
\hskip -1.5cmX\mapsto Hom(X,K) & (\text{на объектах})\\
(X\stackrel {f}{\to }Y)\mapsto (Hom(Y,K)\stackrel {f_*}{\to }Hom(X,K))& (\text{на стрелках})
\end{array}
\right.
$ \\
\vskip 0.05cm
где $f_*(g):=g\circ f$.
\item  {\bf Множество-степень} функтор ${\cal P}:\mathbf {Set}\to \mathbf {Set}:$\\
\vskip 0.05cm
$\left\{
\begin{array}{cr}
\hskip -0.7cmX\mapsto {\cal P}(X) & (\text{на объектах})\\
(X\stackrel {f}{\to }Y)\mapsto ({\cal P}(Y)\stackrel {{\cal P}(f)}{\to }{\cal P}(X))& (\text{на стрелках})
\end{array}
\right.
$ \\
\vskip 0.05cm
где ${\cal P}(X):=\{S\, |\, S\subset X\}$ (множество всех подмножеств),

\addvspace{0.1cm}
${\cal P}(f)(S'):=\{x\in X\, |\, f(x)\in S'\}$ (прообраз подмножества $S'\subset Y$).
\end{enumerate} 

\begin{assignment}
\begin{enumerate}
\item Доказать утверждение 3.
\item Показать, что существует забывающий функтор из категории $\mathbf {Cat}$ в категорию ориентированных графов $\mathbf {Graph}$ (отображение между графами это функция, переводящая вершины в вершины, а ребра в ребра так, что операции взятия начала и конца ребра сохраняются).
\item Убедиться, что произведение функторов одинаковой вариантности дает ковариантный функтор, а разной вариантности -- контравариантный.
\item Найти ковариантную версию функтора множество-степень.
\end{enumerate}
\end{assignment}

\section{Естественные преобразования}

Естественное преобразование можно представлять себе как деформацию одного функтора в другой вдоль стрелок. Содержательно эти стрелки представляют обычно какое-то инвариантное одинаковое для всех объектов правило.

\begin{definition}
Пусть $F:{\cal K}\to {\cal L}$, $G:{\cal K}\to {\cal L}$ два функтора между категориями ${\cal K}$ и ${\cal L}$. Тогда семейство ${\cal L}$-стрелок $\{FX\stackrel {\mu _X}{\to }GX\, |\, X\in {\cal K}_0\}$, по одной для каждого объекта $X$ из ${\cal K}$, называется {\bf естественным преобразованием} функтора $F$ в функтор $G$, если для любой стрелки $(f:X\to Y)\in {\cal K}_1$ диаграмма 
$\vcenter{
\xymatrix{
FY \ar[r]^-{\mu _Y} & GY \\
FX \ar[r]_-{\mu _X} \ar[u]^-{Ff} & GX \ar[u]_-{Gf}
}
}$
коммутативна, то есть $\mu _Y\circ Ff=Gf\circ \mu _X$.
\end{definition}

\noindent Естественные преобразования обозначаются обычно двойными стрелками $\mu :F\Rightarrow G$. Стрелки $\mu _X$ называются {\bf компонентами} естественного преобразования $\mu $.

\addvspace{0.2cm}
{\bf Примеры.}
\begin{enumerate}
\item {\bf Единичное} естественное преобразование функтора $F:{\cal K}\to {\cal L}$ в себя задается единичными стрелками $1_{FX},\ X\in {\cal K}_0$.
\item Рассмотрим два функтора: {\it тождественный} $1_{\mathbf {Ab}}:\mathbf {Ab}\to \mathbf {Ab}$ и {\it функтор кручения} $\mathbf {Tor}:\mathbf {Ab}\to \mathbf {Ab}$. Функтор кручения назначает каждой абелевой группе $A$ подгруппу кручения $\mathbf {Tor}(A):=\{a\in A\, |\, \exists n\in {\Bbb N}, n\ne 0, na=0\}$, а каждой стрелке $A\stackrel {f}{\to }B$ ее ограничение на подгруппу элементов ненулевого кручения $\mathbf {Tor}(A)\stackrel {f|}{\to }\mathbf {Tor}(B)$. Тогда естественным преобразованием является включение групп $\mathbf {Tor}(A)\hookrightarrow A$ для каждой абелевой группы $A$.
\item Пусть $\mathbf {Rng}$ обозначает категорию коммутативных колец и гомоморфизмов. Коммутативное кольцо это множество с двумя коммутативными операциями, сложения и умножения, согласованными между собой дистрибутивным законом $a\cdot (b+c)=a\cdot b+a\cdot c$. Пусть $GL_R(n)$ и $R^*$ обозначают соответственно группу по умножению $n\times n$-матриц с элементами из коммутативного кольца $R$ и группу единиц (обратимых элементов) кольца $R$.\\
Рассмотрим два функтора $Gl_{(-)}(n), (-)^*:\mathbf {Rng}\to \mathbf {Grp}$. Первый сопоставляет кольцу $R$ группу $Gl_R(n)$, а кольцевому гомоморфизму $R\stackrel {f}{\to }S$ гомоморфизм групп $GL_R(n)\to GL_S(n):(a_{ij})\mapsto (f(a_{ij}))$. Второй сопоставляет кольцу $R$ группу $R^*$, а кольцевому гомоморфизму его ограничение на обратимые элементы кольца. Тогда операция взятия определителя матрицы $det:GL_R(n)\to R^*$ является естественным гомоморфизмом групп.
\item Если $Hom(K,-),Hom(K',-):{\cal K}\to \mathbf {Set}$ два ковариантных $Hom$-функтора, то любая стрелка $f:K\to K'$ задает естественное пребразование $Hom(f,-):Hom(K',-)\Rightarrow Hom(K,-)$ по правилу $\forall X\in {\cal K}$ $Hom(f,X):Hom(K',X)\to Hom(K,X):g\mapsto g\circ f$.\\
Аналогично задается естественное преобразование контравариантных $Hom$-функторов $Hom(-,f):Hom(-,K)\to Hom(-,K'):*\mapsto f\circ *$. 
\end{enumerate}

\begin{proposition}
Для любых двух категорий ${\cal K}$ и ${\cal L}$ функторы из ${\cal K}$ в ${\cal L}$ образуют категорию, объектами которой являются сами функторы, а стрелками естественные преобразования между ними с умножением определенным следующим образом:
если $\xymatrix{
{\cal K} \ar@/^1.2pc/[rr]^(0.5955){F}_(0.4){\mu \Downarrow} \ar[rr]^(0.59){G}_(0.42){\nu \Downarrow } \ar@/_1.2pc/[rr]_(0.59){H}           && {\cal L}
}$
два естественных преобразования из $F$ в $G$ и из $G$ в $H$, тогда покомпонентная композиция стрелок составляющих естественные преобразования задает новое естественное преобразование $\nu \circ \mu:F\Rightarrow H$ с компонентами $\nu _X\circ \mu _X$, $X\in {\cal K}^0$. Единицами в категории служат единичные естественные преобразования. 
\end{proposition}
Доказательство.\\
Аксиома единицы очевидно удовлетворяется, как и аксиома ассоциативности, потому что композиция определена покомпонентно.
Проверим является ли композиция $\nu \circ \mu $ естественным преобразованием. Для этого внешний прямоугольник должен быть коммутативен

\addvspace{0.1cm}
$\vcenter{\xymatrix{
FY \ar[r]^-{\mu _Y}  & GY \ar[r]^-{\nu _Y} & HY\\
FX \ar[u]^-{Ff} \ar[r]_-{\mu _X}  & GX \ar[u]^-{Gf} \ar[r]_-{\nu _X}  & HX \ar[u]_-{Hf}
}}$ что также выполняется, потому что составляющие его квадраты коммутативны по условию. \hfill  $\square $

\addvspace{0.2cm}
Категрия функторов из ${\cal K}$ в ${\cal L}$ обозначается через степень ${\cal L}^{\cal K}$. Контравариантные функторы также составляют (свою) категорию ${\cal L}^{{\cal K}^{op}}$.

\begin{proposition}
Для любой категории ${\cal K}$ имеются 
\begin{itemize}
\item контравариантный функтор ${\cal K}^{op}\to {\mathbf {Set}}^{\cal K}$,
\item ковариантный функтор ${\cal K}\to {\mathbf {Set}}^{{\cal K}^{op}}$,

\addvspace{0.1cm}
где ${\mathbf {Set}}^{\cal K}$, ${\mathbf {Set}}^{{\cal K}^{op}}$ соответственно категории ковариантных и контравариантных функторов из ${\cal K}$ в $\mathbf {Set}$. 
\end{itemize}
\end{proposition}
Доказательство.\\
Назначаем объектам и стрелкам $Hom$-функторы как в примере 4. Проверка того, что это дает функторы, остается в качестве упражнения.              \hfill $\square $

\begin{assignment}
\begin{enumerate}
\item Убедиться, что в приведенных примерах действительно определены естественные преобразования.
\item Нарисовать коммутативный квадрат для естественного преобразования контравариантных функторов. 
\item Восполнить недостающие детали в доказательствах утверждений.
\end{enumerate}
\end{assignment}

\section{2-категории}

Естественней и проще рассматривать категорию $\mathbf {Cat}$ как состоящую из категорий в качестве объектов, функторов в качестве стрелок и естественных преобразований в качестве двумерных стрелок (2-стрелок). Это приводит к понятию 2-категории.

\begin{definition}
2-категория ${\cal K}$ состоит из трех непересекающихся классов ${\cal K}_0$ (объекты), ${\cal K}_1$ (1-стрелки), ${\cal K}_2$ (2-стрелки), при этом
\begin{itemize}
\item для элементов классов ${\cal K}_0$ и ${\cal K}_1$ существуют единичные стрелки, на единицу большей размерности, чем исходный объект или стрелка,   
\item на классах ${\cal K}_1$ и ${\cal K}_2$ определены операции $d$ (взятия начала) и $c$ (взятия конца) стрелки, такие что $dc=dd$, $cd=cc$ (смотрите картинку $\xymatrix {\bullet \ar@/^0.5pc/[r]_-{\Downarrow } \ar@/_0.5pc/[r] & \bullet }$),
\item на 1-стрелках и 2-стрелках определено обычное умножение $f\circ g$ (или $\mu \circ \nu$), если конец первой стрелки совпадает с началом второй, то есть $cg=df$ (или $c\nu =d\mu $) (на диаграммах это будет $\xymatrix{\bullet \ar[r]^-{g} & \bullet \ar[r]^-{f} & \bullet }$ или $\xymatrix{\bullet \ar@/^1pc/[r]_-{\nu \Downarrow } \ar[r]_-{\mu \Downarrow } \ar@/_1pc/[r]   &  \bullet }$). При этом умножение $\circ $ для 1-стрелок называется {\bf горизонтальным}, а для 2-стрелок {\bf вертикальным},
\item для 2-стрелок определено {\bf горизонтальное} умножение $\mu *\nu $, если $dd\mu =cc\nu $ (то есть $\xymatrix{\bullet  \ar@/^0.5pc/[r]^-{F}_-{\nu \Downarrow } \ar@/_0.5pc/[r]_-{G} & \bullet \ar@/^0.5pc/[r]^-{H}_-{\mu \Downarrow } \ar@/_0.5pc/[r]_-{J}  & \bullet }$), при этом
\begin{itemize}
\item тип произведения есть $\mu *\nu :H\circ F\Rightarrow J\circ G$,
\item $1_H*1_F=1_{H\circ F}$ для любых умножаемых 1-стрелок $F$ и $H$, 
\end{itemize}
\item (аксиома единицы) для $f:X\to Y$ и для $\mu :f\Rightarrow g$ выполняется $1_Y\circ f=f=f\circ 1_X$, $1_g\circ \mu =\mu =\mu \circ 1_f$, $1_{1_Y}*\mu =\mu =\mu *1_{1_X}$,
\item (аксиома ассоциативности) $f\circ (g\circ h)=(f\circ g)\circ h$, $\mu \circ (\nu \circ \eta )=(\mu \circ \nu )\circ \eta$, $\alpha *(\beta *\gamma )=(\alpha *\beta )*\gamma $ (для всех стрелок, для которых произведения определены),
\item (перестановка операций) $(\delta \circ \gamma )*(\beta \circ \alpha )=(\delta *\beta )\circ (\gamma *\alpha )$ для 2-стрелок, расположенных следующим образом $\xymatrix{\bullet \ar@/^1pc/[r]_-{\alpha  \Downarrow } \ar[r]_-{\beta  \Downarrow } \ar@/_1pc/[r]  & \bullet \ar@/^1pc/[r]_-{\gamma \Downarrow } \ar[r]_-{\delta \Downarrow } \ar@/_1pc/[r]  & \bullet }$
\end{itemize}
\end{definition}
Это длинное определение вводит естественные правила, с помощью которых можно производить вычисления со стрелками. По существу, утверждается, что вычисление с одной операцией $\circ $ или $*$ не зависит от расстановки скобок, что единица поглощается обычным образом, и что операции $\circ $ и $*$ коммутируют (перестановочны). Другое определение можно найти в упражнении 1.

\addvspace{0.2cm}
{\bf Примеры.}
\begin{enumerate}
\item Категория $\mathbf {Cat}$ является 2-категорией с горизонтальным умножением для 2-стрелок $\xymatrix{\bullet \ar@/^0.5pc/[r]^-{F}_-{\alpha \Downarrow } \ar@/_0.5pc/[r]_-{G} & \bullet \ar@/^0.5pc/[r]^-{H}_-{\beta \Downarrow } \ar@/_0.5pc/[r]_-{J} & \bullet }$, определенным следующим образом $(\beta *\alpha )_X:=J(\alpha _X)\circ \beta _{FX}$ (что также равно $\beta _{GX}\circ H(\alpha _X)$, поскольку $\beta $ естественное преобразование. Действительно, квадрат $\vcenter{\xymatrix{HFX \ar[r]^-{\beta _{FX}} \ar[d]_-{H\alpha _X} & JFX \ar[d]^-{J\alpha _X}\\
HGX \ar[r]_-{\beta _{GX}} & JGX }}$ коммутативный по определению естественного преобразования $\beta $).\\
Доказательство ассоциативности $*$ и перестановочного закона для $*$ и $\circ $ оставляем читателю в качестве упражнения.

$\mathbf {Cat}$, рассматриваемая как 2-категория, обозначается $2\text{-}\mathbf {Cat}$.
\item $Hom$-функтор $2\text{-}\mathbf {Cat}({\cal K},-):2\text{-}\mathbf {Cat}\to 2\text{-}\mathbf {Cat}$ является так называемым 2-функтором (который отображает объекты в объекты, 1-стрелки в 1-стрелки, 2-стрелки в 2-стрелки, сохраняет единицы и обе операции). Он назначает категории ${\cal L}$ 1-категорию функторов из ${\cal K}$ в ${\cal L}$, то есть $2\text{-}\mathbf {Cat}({\cal K},{\cal L})$, функтору $F:{\cal L}\to {\cal M}$ функтор $2\text{-}\mathbf {Cat}({\cal K},F):2\text{-}\mathbf {Cat}({\cal K},{\cal L})\to 2\text{-}\mathbf {Cat}({\cal K},{\cal M}):G\mapsto F\circ G$, и естественному преобразованию $\mu :F\Rightarrow H$, где $F,H:{\cal K}\to {\cal L}$ функторы, естественное преобразование $2\text{-}\mathbf {Cat}({\cal K},\mu ):G\mapsto \mu *1_G$, где $G$ пробегает объекты (функторы) в $2\text{-}\mathbf {Cat}({\cal K},{\cal L})$.  
\item $2\text{-}\mathbf {Top}$ состоит из топологических пространств в качестве объектов, непрерывных отображений в качестве 1-стрелок и гомотопий (точнее, классов эквивалентности гомотопий) в качестве 2-стрелок. Гомотопия (=деформация) $\alpha :f\Rightarrow g$ между двумя непрерывными отображениями $f:X\to Y$ и $g:X\to Y$ определяется как 'путь' $\alpha _t:X\to Y, \ t\in [0,1]$, такой что $\alpha _0(x)=f(x)$, $\alpha _1(x)=g(x)$ для всех точек $x\in X$, и функция $\alpha _-(-):[0,1]\times X\to Y$ непрерывна.\\
Горизонтальное умножение для непрерывных функций обычная композиция $(f\circ g):=f(g(x)),\ x\in X$, вертикальное умножение для гомотопий $(\beta \circ \alpha )_t(x):=\left\{ 
\begin{array}{lr}
\alpha _{2t}(x) & t\in [0,1/2]\\
\beta _{2t-1}(x) & t\in [1/2,1]
\end{array}  
\right.$
(то есть первую половину 'времени' выполняется $\alpha $, а вторую $\beta $).
Горизонтальное умножение $*$ для гомотопий похоже на соответствующее умножение в $2\text{-}\mathbf {Cat}$ и оставляется в качестве упражнения. 
\item 2-функтор $2\text{-}\mathbf {Top}(\mathbf {1},-):2\text{-}\mathbf {Top}\to 2\text{-}\mathbf {Cat}$, где $\mathbf 1$ одноточечное топологическое пространство, называется функтором взятия {\bf фундаментального группоида}. Он назначает топологическому пространству 2-категорию, объектами которой являются точки пространства, 1-стрелками непрерывные пути между точками, 2-стрел-\\ ками гомотопии (классы эквивалентности гомотопий) между путями. Этот функтор имеет важные применения в топологии. 
\end{enumerate}

\begin{assignment}
\begin{enumerate}
\item Существует более абстрактное определение 2-категории. Именно, 2-категория ${\cal K}$ состоит из класса объектов ${\cal K}_0$, для каждых двух элементов которого $X,Y\in {\cal K}_0$ назначена категория ${\cal K}(X,Y)$, объекты которой называются 1-стрелками, а стрелки 2-стрелками. Для каждой категории вида ${\cal K}(X,X)$ имеется функтор $I:\mathbf {1}\to {\cal K}(X,X)$, где $\mathbf 1$ категория состоящая из одного объекта и одной стрелки, (выделяющий $1_X$ и $1_{1_X}$ в ${\cal K}(X,X)$). Для каждой тройки объектов $X,Y,Z$ имеется функтор умножения $*:{\cal K}(Y,Z)\times {\cal K}(X,Y)\to {\cal K}(X,Z)$. При этом выполняются 

(аксиома единицы) 

$\vcenter{\xymatrix{
{\cal K}(X,X)\times {\cal K}(X,Y) \ar[r]^-{*} & {\cal K}(X,Y)\\
\mathbf 1\times {\cal K}(X,Y) \ar[u]^-{I\times 1_{{\cal K}(X,Y)}} \ar[ur]_-{pr_2} &
}}$ 
$\vcenter{\xymatrix{
{\cal K}(X,Y)\times {\cal K}(X,X) \ar[r]^-{*} & {\cal K}(X,Y)\\
\mathbf {\cal K}(X,Y)\times \mathbf 1 \ar[u]^-{1_{{\cal K}(X,Y)}\times I} \ar[ur]_-{pr_1} &
}}$ 

и (аксиома ассоциативности)

$\xymatrix{
{\cal K}(Z,W)\times ({\cal K}(Y,Z)\times {\cal K}(X,Y)) \ar[r]^-{\stackrel {\sim }{\to }} \ar[d]_-{1_{{\cal K}(Z,W)}\times *} & ({\cal K}(Z,W)\times {\cal K}(Y,Z))\times {\cal K}(X,Y) \ar[d]_-{*\times 1_{{\cal K}(X,Y)}}\\
{\cal K}(Z,W)\times {\cal K}(X,Z) \ar[d]_-{*} & {\cal K}(Y,W)\times {\cal K}(X,Y) \ar[dl]^-{*}\\
{\cal K}(X,W) & 
}$

В определении используются произведения категорий вида ${\cal K}\times {\cal L}$ и произведения функторов вида $F\times G$. Не давая формального определения, только отметим, что объекты произведения категорий это пары объектов соответствующих категорий, стрелки пары стрелок, и все операции выполняются покомпонентно. Произведение функторов действует по формуле $(F\times G)(K,L):=(FK,GL)$ (на объектах) и $(F\times G)(f,g):=(Ff,Gg)$ (на стрелках). Функторы проекции $pr_1:{\cal K}\times {\cal L}\to {\cal K}$ и $pr_2:{\cal K}\times {\cal L}\to {\cal L}$ выделяют соответственно первую и вторую компоненту пары (объектов или стрелок).

С точки зрения этого определения горизонтальное умножение 1-стрелок следовало бы обозначать через $*$.

Читателю предлагается проанализировать это определение и убедиться, что оно эквивалентно данному в тексте.
\item Доказать ассоциативность и перестановочный закон для $*$ в категории $2\text{-}\mathbf {Cat}$.
\item Выписать явную формулу для горизонтального умножения гомотопий.
\item Убедиться, что функтор взятия фундаментального группоида назначает непрерывным отображениям функторы, а гомотопиям естественные преобразования.
\end{enumerate}
\end{assignment}

\addvspace{0.1cm}
{\bf Замечание.} Обычные категории могут быть наделены структурой 2-категории по крайней мере двумя разными способами. Когда вводятся только единичные 2-стрелки для каждой 1-стрелки категории. И когда считается, что между любыми двумя параллельными стрелками (то есть с общим началом и концом) существует одна и только одна 2-стрелка.

\section{Эквивалентности в 2-категории}

Часто в математике интересуются вопросом, существует ли {\it обратимое} преобразование одного объекта в другой, при этом преобразования могут быть более или менее обратимыми. Соответственно вводятся строгие эквивалентности (или изоморфизмы) и слабые эквивалентности.

\begin{definition}
1-стрелка $f:X\to Y$ в 2-категории ${\cal K}$ называется
\begin{itemize}
\item {\bf строгой эквивалентностью}, если существует обратная стрелка $g:Y\to X$, то есть такая стрелка, что $g\circ f=1_X$, $f\circ g=1_Y$, 
\item {\bf слабой эквивалентностью}, если существует стрелка $g:Y\to X$, и существуют обратимые 2-стрелки $\alpha $ и $\omega $ такие, что $g\circ f\stackrel {\alpha }{\Rightarrow }1_X$, $f\circ g\stackrel {\omega }{\Rightarrow }1_Y$.
\end{itemize}
\end{definition}

\noindent В определении существенно, что $\alpha $ и $\omega $ обратимы (то есть 2-изоморфизмы).

\noindent Сильные эквивалентности обозначаются $X\stackrel {\simeq }{\to }Y$, а слабые $X\stackrel {\sim }{\to }Y$. Соответственно изоморфные объекты обозначаются $X\simeq Y$, а (слабо) эквивалентные $X\sim Y$.

Очевидно, строгие эквивалентности являются слабыми (достаточно положить $\alpha =1_{1_X}$, $\omega =1_{1_Y}$), но не наоборот.

\addvspace{0.2cm}
{\bf Примеры.}
\begin{enumerate}
\item В категории множеств $\mathbf {Set}$, рассматриваемой как 2-категория с единичными только 2-стрелками, все строгие эквивалентности совпадают со слабыми и являются биекциями множеств. Если рассматривать $\mathbf {Set}$ с 2-стрелками, существующими по одной для каждой пары параллельных стрелок, то строгими эквивалентностями будут биекции, а слабыми любые отображения между непустыми множествами. В последнем случае все непустые множества эквивалентны.
\item В категории $2\text{-}\mathbf {Cat}$ отношение изоморфизма является слишком жестким, но тем не менее полезным. Оно выражает факт, что две категории 'абсолютно' одинаковы, хотя, может быть, описаны по разному. Так категория $\mathbf {Ab}$ изоморфна категории $\Bbb Z$-модулей $\Bbb Z\text{-}\mathbf {Mod}$. $\Bbb Z$-модуль это абелева группа, элементы которой можно умножать на целые чила, с простыми аксиомами. Казалось бы, дополнительная структура умножения на числа должна делать класс $\Bbb Z$-модулей меньше, чем класс абелевых групп, но это умножение выводится из внутренней структуры группы, так что каждая абелва группа есть $\Bbb Z$-модуль ($n\cdot a=a+\cdots +a$ ($n$ раз)).\\
Отношение слабой эквивалентности в $2\text{-}\mathbf {Cat}$ называется просто отношением эквивалентности категорий. Согласно определению, если $\xymatrix{{\cal K} \ar@/^0.5pc/[r]^-{F} & {\cal L} \ar@/^0.5pc/[l]^-{G}}$ эквивалентность категорий ${\cal K}$ и ${\cal L}$, то $X\simeq GFX$, $Y\simeq FGY$ для всех объектов $X\in {\cal K}_0$, $Y\in {\cal L}_0$.  Например, категория векторных пространств $\mathbf {Vect}$ эквивалентна своей полной подкатегории, состоящей из канонических представителей векторных пространств для каждой размерности $0,\ k,\ k^2,\dots ,k^n,\dots $, где $k$ есть {\it поле} чисел, на которые умножаются элементы векторных пространств.
\item Возьмем два объекта в $2\text{-}\mathbf {Top}$, именно, прямую линию и точку. Тогда эти два объекта не изоморфны в $2\text{-}\mathbf {Top}$, так как изоморфные объекты в $2\text{-}\mathbf {Top}$ должны быть изоморфны в $2\text{-}\mathbf {Set}$ (с 2-стрелками, выбранными по одной для каждой пары параллельных стрелок в $\mathbf {Set}$) на том основании, что мы применяем забывающий функтор. Но точки прямой и одна точка не могут быть во взаимно-однозначном соответствии. Тем не менее, прямая и точка слабо эквивалентны в $2\text{-}\mathbf {Top}$. Доказательство остается читателю.
\end{enumerate}

\begin{proposition}
Отношения изоморфизма и слабой эквивалентности являются отношениями эквивалентности (то есть рефлексивны, симметричны и транзитивны).
\end{proposition}
Доказательство.
Для изоморфизма это уже доказано в упражнении 2.5 параграфа 1.2. Докажем предложение для слабой эквивалентности.\\
(рефлексивность) $X\sim X$, для этого выбираем 1-стрелки $\xymatrix{X \ar@/^0.5pc/[r]^-{1_X} & X \ar@/^0.5pc/[l]^-{1_X}}$ и 2-стрелки $\alpha =\beta =1_{1_X}$,

\addvspace{0.2cm}
\noindent (симметричность) $X\sim Y$ $\Rightarrow $ $Y\sim X$, для этого меняем местами $\alpha $ и $\beta $,

\addvspace{0.2cm}
\noindent (транзитивность) $X\sim Y$ и $Y\sim Z$ дают $X\sim Z$, для этого рассмотрим диаграмму
$\xymatrix{X \ar@/^0.5pc/[r]^-{F}  & Y \ar@/^0.5pc/[r]^-{H}  \ar@/^0.5pc/[l]^-{G} & Z \ar@/^0.5pc/[l]^-{J} }$ и вычислим
$1_X\stackrel {\alpha _X^{-1}}{\Rightarrow } G\circ F\stackrel {=}{\Rightarrow }G\circ 1_Y\circ F\stackrel {1_G*\alpha _Y^{-1}*1_F}{\Longrightarrow }G\circ J\circ H\circ F$, то есть $1_X\simeq (G\circ J)\circ (H\circ F)$. Аналогично,
$1_Z\simeq (H\circ F)\circ (G\circ J)$. Следовательно, $X\sim Z$.          \hfill $\square $

\begin{assignment}
\begin{enumerate}
\item Доказать все утверждения примера 1.
\item Доказать, что прямая и точка слабо эквивалентны в категории $2\text{-}\mathbf {Top}$.
\item Доказать, что вертикальная и горизонтальная композиции 2-изоморфизмов дают 2-изоморфизмы.
\item Доказать, что для 1-стрелок $F:Y\to Z$ и $G:W\to X$ горизонтальные умножения \\
$F*-:{\cal K}(X,Y)\to {\cal K}(X,Z):
\left\{
\begin{array}{cr}
H\mapsto F\circ H & (\text{на объектах})\\
\mu \mapsto 1_F*\mu & (\text{на стрелках})
\end{array}
\right.
$ 
\ \ и \\
$-*G:{\cal K}(X,Y)\to {\cal K}(W,Y):
\left\{
\begin{array}{cr}
H\mapsto H\circ G & (\text{на объектах})\\
\mu \mapsto \mu *1_G & (\text{на стрелках})
\end{array}
\right.
$ 
являются функторами.
\item Доказать, что вложение полной подкатегории ${\cal L}\hookrightarrow {\cal K}$ такой, что каждый объект из ${\cal K}$ изоморфен некоторому объекту в ${\cal L}$, является эквивалентностью.
\end{enumerate}
\end{assignment}

\chapter{Представимые функторы}

Представимые функторы являются самым старомодным средством в теории категорий, но удобны тем, что позволяют рассматривать с единой точки зрения разные вопросы (пределы, копределы, сопряженность и другие).

\section{Представимые функторы. Вложение Йонеды}

\begin{definition}
\begin{itemize}
\item Ковариантный функтор $F:{\cal K}\to \mathbf {Set}$ называется {\bf представимым}, если он изоморфен $Hom$-функтору \,  ${\cal K}(A,-):{\cal K}\to \mathbf {Set}$ для некоторого объекта $A\in {\cal K}_0$, то есть если существует естественное проебразование $\chi :{\cal K}(A,-)\Rightarrow F$ все компоненты которого $\chi _X:{\cal K}(A,X)\stackrel {\simeq }{\to }F(X)$ биекции. 
\item Аналогично, контравариантный функтор $F:{\cal K}^{op}\to \mathbf {Set}$ называется {\bf представимым}, если он изоморфен $Hom$-функтору \, ${\cal K}(-,A):{\cal K}^{op}\to \mathbf {Set}$ для некоторого объекта $A\in {\cal K}_0$, то есть если существует естественное проебразование $\chi :{\cal K}(-,A)\Rightarrow F$ все компоненты которого $\chi _X:{\cal K}(X,A)\stackrel {\simeq }{\to }F(X)$ биекции. 
\end{itemize}
\end{definition}

\noindent Объект $A$ называется {\bf представляющим объектом}. Множество $FA$ содержит характерный элемент $\chi _A(1_A)$, который называется {\bf универсальным элементом}.

Из определения универсального элемента и того, что $\chi $ естественное преобразование, следует, что все компоненты $\chi $ полностью определены самим функтором $F$ и универсальным элементом $\chi _A(1_A)$. Действительно, рассмотрим диаграмму (например, для контравариантного функтора $F$)
\vskip 0.15cm
\noindent $\vcenter{\xymatrix{
X \ar[d]_-{f}  & {\cal K}(X,A)  \ar[r]^-{\chi _X}      &      FX  \\  
A              & {\cal K}(A,A)  \ar[r]_-{\chi _A} \ar[u]^-{{\cal K}(f,A)}     &      FA  \ar[u]_-{F(f)}
}}$ 
\hskip 0.5cmПолучим $\chi _X\circ {\cal K}(f,A)(1_A)=Ff\circ \chi _A(1_A)$ 
\vskip 0.15cm
\noindent или $\chi _X(f)=F(f)(\chi _A(1_A))$.
Из последней формулы также следует, что функтор $F$ со значениями в категории множеств представим тогда и только тогда, когда существуют объект $A\in {\cal K}_0$ и элемент $a\in FA$ такие, что отображение $F(-)(a):{\cal K}(X,A)\to FX$ биективно для каждого $X\in {\cal K}_0$ (для ковариантного функтора отображение будет типа $F(-)(a):{\cal K}(A,X)\to FX$).

Если ввести 'категорию элементов', объектами которой являются пары $(C,c)$, где $C\in {\cal K}_0$, $c\in FC$, а стрелками $f:(B,b)\to (C,c)$ служат стрелки $(f:B\to C)\in {\cal K}_1$ такие, что $F(f)(c)=b$ (для контравариантного функтора $F$) или $F(f)(b)=c$ (для ковариантного функтора $F$), то контравариантный функтор $F:{\cal K}^{op}\to \mathbf {Set}$ представим объектом $A$ с универсальным элементом $a$ тогда и только тогда, когда $(A,a)$ является {\it конечным объектом} в категории элементов, и соответственно ковариантный функтор $F:{\cal K}\to \mathbf {Set}$ представим объектом $A$ с универсальным элементом $a$ тогда и только тогда, когда $(A,a)$ является {\it начальным объектом} в категории элементов.

Функтор $F:{\cal K}\to \mathbf {Set}$ может быть представим разными объектами из $\cal K$ с соответственно разными универсальными элементами. Однако, в каждом из этих случаев представляющая пара является начальным объектом в категории элементов, и это значит, что если $(A,a)$ и $(A',a')$ две представляющие пары для функтора $F$, то существует единственный изоморфизм между ними $f:(A,a)\stackrel {\simeq }{\to }(A',a')$. Аналогично, для контравариантного представимого функтора $F:{\cal K}^{op}\to \mathbf {Set}$ пара $(A,a)$ определена {\it однозначно с точностью до изоморфизма}. 

\addvspace{0.2cm}
{\bf Примеры.}
\begin{enumerate}
\item $Hom$-функтор ${\cal K}(A,-)$ является представимым с представляющим объектом $A$, так как он изоморфен самому себе. Вложение категории $Hom$-функторов в категорию всех представимых функторов является эквивалентностью.  
\item Контравариантный функтор множество-степень ${\cal P}:\mathbf {Set}\to \mathbf {Set}$ представим двухэлементным множеством, например, $\{0,1\}$ (в то же время ковариантный функтор множество-степень не является представимым). Универсальным элементом является одноэлементное подмножество $\{0\}\in {\cal P}(\{0,1\})$ или $\{1\}\in {\cal P}(\{0,1\})$. Действительно, ${\cal P}(-)(\{0\}):\mathbf {Set}(X,\{0,1\})\to {\cal P}(X):f\mapsto f^{-1}(\{0\})$ является биекцией (обратное отображение сопоставляет подмножеству множества $X$ функцию, равную нулю на этом подмножестве и единице вне его).
\item Забывающий функтор $U:\mathbf {Top}\to \mathbf {Set}$ является представимым одноточечным топологическим пространством $\mathbf {1}$. Для универсального элемента нет выбора, кроме самой этой точки. Обозначим ее $*$ ($*\in U(\mathbf 1)$). Проверяем, $U(-)(*):\mathbf {Top}(\mathbf 1,X)\to U(X):f\mapsto f(*)$ является биекцией (обратное преобразование назначает элементу $x\in U(X)\in \mathbf {Set}_0$ непрерывную функцию $f:\mathbf 1\to X$ с единственным значением $x$).
\item Забывающий функтор $U:\mathbf  {Vect}\to \mathbf {Set}$ является представимым. Представляющий объект есть $k$ (поле скаляров), рассматриваемое как одномерное векторное пространство. Универсальным элементом является единица поля $1\hskip -0.05cm\in \hskip -0.05cmk$. Проверяем, $U(-)(1):\hskip -0.05cm\mathbf {Vect}(k,X)\hskip -0.05cm\to U(X):f\mapsto f(1)$ является биекцией, так как для любого элемента $x$ векторного пространства $X$ найдется единственная линейная функция, отображающая $1\in k$ в $x$.
\end{enumerate}

Задавая формулу естественного преобразования $F(-)(a)$ $Hom$-функ-\\ тора ${\cal K}(A,-)$ в представимый функтор $F$, мы никак не учитывали заранее, что это естественное преобразование будет изоморфизмом (это являлось дополнительным условием для проверки, что объект $A$ и элемент $a$ определяют представимый функтор $F$). Это значит, что произвольное естественное преобразование $Hom$-функтора ${\cal K}(A,-)$ в любой функтор $F$ со значениями в категории множеств выглядит точно также $F(-)(a)$.

\begin{proposition}[Лемма Йонеды]
Существует взаимно-однозначное соответсвие $\theta _{F,A}:{\mathbf {Set}}^{\cal K}({\cal K}(A,-),F)\stackrel {\simeq }{\to }F(A)$ (между естественными преобразованиями любого $Hom$-функтора в произвольный функтор $F$ со значениями в $\mathbf {Set}$ и множеством $F(A)$, которое функтор $F$ назначает представляющему объекту $A$). Это соответствие естественное по обоим аргументам $F\in (\mathbf {Set}^{\cal K})_0$  и $A\in {\cal K}_0$. 
\end{proposition}
Доказательство.\\
Определим $\theta _{F,A}$ тем же способом, каким определяли универсальный элемент, то есть $\theta _{F,A}(\alpha ):=\alpha _A(1_A)$, а обратное преобразование формулой $\bar \theta _{F,A}(a):=F(-)(a)$. Докажем, что $\theta _{F,A}$ и $\bar \theta _{F,A}$ действительно взаимно-обратные преобразования:\\
$\theta _{F,A}\bar \theta _{F,A}(a)=\theta _{F,A}(F(-)(a))=F(1_A)(a)=1_{FA}(a)=a$,\\
$\bar \theta _{F,A}\theta _{F,A}(\alpha )=\bar \theta _{F,A}(\alpha _A(1_A))=F(-)(\alpha _A(1_A))=\alpha $, то есть $\bar \theta _{F,A}=(\theta _{F,A})^{-1}$.
Докажем естественность $\theta _{F,A}$ только по аргументу $A$. Рассмотрим диа-\\

\addvspace{0.0cm}
\noindent грамму \hskip 1.0cm
$\vcenter{\xymatrix{
A \ar[d]_-{f} & \mathbf {Set}^{\cal K}({\cal K}(A,-),F) \ar[r]^-{\theta _{F,A}} \ar[d]_-{\mathbf {Set}^{\cal K}({\cal K}(f,-),F)} &  F(A) \ar[d]^-{F(f)}\\
B & \mathbf {Set}^{\cal K}({\cal K}(B,-),F) \ar[r]_-{\theta _{F,B}} &  F(B) 
}}$\hskip 1.0cmВычисляем 

\addvspace{0.2cm}
\noindent $\theta _{F,B}({\mathbf {Set}^{\cal K}({\cal K}(f,-),F)}(\mu ))=\theta _{F,B}(\mu \circ {\cal K}(f,-))=\mu _B\circ ({\cal K}(f,-))_B\, (1_B)=$

\noindent $\mu _B(({\cal K}(f,-))_B(1_B))=\mu _B(1_B\circ f)=\mu _B(f\circ 1_A)=\mu _B({\cal K}(f,-)(1_A))=\mu _B\circ {\cal K}(f,-)\, (1_A)=[\mu \text{ естественное преобразование }{\cal K}(A,-) \text{ в } F]=F(f)\circ \mu _A\, (1_A)=F(f)(\theta _{F,A}(\mu ))$, то есть диаграмма коммутативна, и $\theta _{F,A}$ естественное преобразование в аргументе $A$.               \hfill $\square $

Лемма Йонеды дает хорошее описание множества всех естественных преобразований из $Hom$-функтора в некоторый отмеченный функтор и действительно полезна, когда такая задача возникает (например, в топологии, в теории характеристических классов).

Прямым приложением леммы является доказательство того, что функторы $X\mapsto {\cal K}(X,-)$ и $X\mapsto {\cal K}(-,X)$ являются биективными, то есть осуществляют вложение (соответственно контравариантное или ковариантное) категории ${\cal K}$ как подкатегории категории функторов со значениями в $\mathbf {Set}$.

\begin{proposition}[Вложение Йонеды]
\begin{itemize}
\item Ковариантный функтор 

$\mathbf {Y}:{\cal K}\to \mathbf {Set}^{{\cal K}^{op}}:$

\addvspace{0.2cm}
$\left\{
\begin{array}{c}
\hskip -4.0cmA\mapsto {\cal K}(-,A)\\
(f:A\to B)\mapsto {\cal K}(-,f):{\cal K}(-,A)\to {\cal K}(-,B):g\mapsto f\circ g 
\end{array}
\right. $ 

\addvspace{0.2cm}
является инъективным на объектах и биективным.
\item Аналогично, контравариантный функтор 

$\bar {\mathbf {Y}}:{\cal K}\to \mathbf {Set}^{{\cal K}}:$

\addvspace{0.2cm}
$\left\{
\begin{array}{c}
\hskip -4.0cmA\mapsto {\cal K}(A,-)\\
(f:A\to B)\mapsto {\cal K}(f,-):{\cal K}(B,-)\to {\cal K}(A,-):g\mapsto g\circ f 
\end{array}
\right. $ 

\addvspace{0.2cm}
является инъективным на объектах и биективным.
\end{itemize}
\end{proposition}
Доказательство.\\
Рассмотрим для определенности только ковариантный случай. Инъективность на объектах следует из того, что если $A\ne B$, то для любого объекта $X\in {\cal K}_0$ множества ${\cal K}(X,A)$ и ${\cal K}(X,B)$ различны, так как состоят из стрелок с разными концами (хотя эти множества могут быть во взаимно-однозначном соответствии).\\
Подставим в лемму Йонеды $F={\cal K}(-,B)$, получим\\ 
$\theta _{{\cal K}(-,B),A}:\mathbf {Set}^{{\cal K}^{op}}({\cal K}(-,A),{\cal K}(-,B))\stackrel {\simeq }{\to }{\cal K}(A,B)$, и как было определено раньше 
$\theta ^{-1}_{{\cal K}(-,B),A}(f):={\cal K}(-,B)(f)$. Напишем формулу для $X\text{-}$компоненты этого естественного преобразования\\
$(\theta ^{-1}_{{\cal K}(-,B),A}(f))_X:=({\cal K}(-,B)(f))_X:{\cal K}(X,A)\to {\cal K}(X,B):g\mapsto f\circ g=(\mathbf Y(f))_X(g)$, где $f\in {\cal K}(A,B)$, $g\in {\cal K}(X,A)$, то есть $(\theta ^{-1}_{{\cal K}(-,B),A}(f))_X=(\mathbf Y(f))_X$ для каждого $X\in {\cal K}_0$. Это значит, что $\theta ^{-1}_{{\cal K}(-,B),A}(f)=\mathbf Y(f)$ для каждого $f\in {\cal K}(A,B)$, то есть отображение Йонеды биективно на $Hom$-множествах.    \hfill  $\square $

\begin{assignment}
\begin{enumerate}
\item Показать, что пара $(\{0,1\},\{0,1\})$, где $\{0,1\}\in {\cal P}(\{0,1\})=\{\O , \{0\}, \{1\}, \{0,1\}\}$, не является представляющей для функтора множества-степень ${\cal P}$.
\item Найти представляющий объект и универсальный элемент забывающего функтора $U:\mathbf {Grp}\to \mathbf {Set}$.
\item Доказать, что формула $F(-)(a)$ действительно определяет естественное преобразование ${\cal K}(A,-)\to F$, $a\in F(A)$ (для ковариантного функтора $F$) или ${\cal K}(-,A)\to F$, $a\in F(A)$ (для контравариантного функтора $F$).  
\end{enumerate}
\end{assignment}

\section{Пределы}

Пределы представляют собой объекты, универсально присоединенные к диаграммам. Диаграмма это просто набор объектов и некоторых стрелок между ними, которые необязательно составляют категорию. Говоря о пределах (и копределах), мы не будем предполагать, что диаграммы коммутативны. Удобно задавать диаграмму отображением подходящего ориентированного графа в категорию.Требуется только, чтобы вершины переходили в вершины, ребра в стрелки, и чтобы начало и конец ребра переходили в начала и конец стрелки-образа. Например, для графа, состоящего только из одной вершины \, $\bullet $\, , диаграмма это просто один объект, например, \, $X$. Для графа, состоящего из двух вершин \, $\bullet \hskip 0.5cm \bullet$\, , диаграммой будет пара объектов, например, \, $X \hskip 0.5cmY$. Для графа $\xymatrix {\bullet \ar@/^0.5pc/[r] \ar@/_0.5pc/[r] & \bullet }$  диаграммой будет, например, $\xymatrix {X \ar@/^0.5pc/[r]^-{f} \ar@/_0.5pc/[r]_-{g} & Y}$, и так далее. Естественное преобразование диаграмм определяется также как для функторов набором стрелок, по одной для каждой вершины графа так, что эти стрелки коммутируют со стрелками диаграммы. Например, для последнего случая естественное преобразование диаграммы $\xymatrix {X \ar@/^0.5pc/[r]^-{f} \ar@/_0.5pc/[r]_-{g} & Y}$ в диаграмму $\xymatrix {X' \ar@/^0.5pc/[r]^-{f'} \ar@/_0.5pc/[r]_-{g'} & Y'}$  задается парой стрелок $\alpha _1:X\to X'$ и $\alpha _2:Y\to Y'$ таких, что квадраты
$\vcenter{\xymatrix{
X  \ar[r]^-{\alpha _1} \ar[d]_-{f} &   X' \ar[d]^-{f'} \\
Y  \ar[r]_-{\alpha _2} &   Y'
}}$ 
и
$\vcenter{\xymatrix{
X  \ar[r]^-{\alpha _1} \ar[d]_-{g} &   X' \ar[d]^-{g'} \\
Y  \ar[r]_-{\alpha _2} &   Y'
}}$
коммутативны.

\begin{proposition}
\begin{itemize}
\item Диаграммы из графа $G$ в категорию ${\cal K}$ составляют категорию, с объектами -- диаграммами и стрелками -- естественными преобразованиями диаграмм. Эта категория обозначается $\mathbf {Dgrm}_{G,{\cal K}}$.
\item Имеется функтор $\Delta :{\cal K}\to \mathbf {Dgrm}_{G,{\cal K}}$, назначающий объекту $X\in {\cal K}_0$ постоянную диаграмму $\Delta (X)$ (отображающую все вершины графа $G$ в объект $X$ и все стрелки графа $G$ в единичную стрелку $1_X$) и стрелке $f:X\to Y$ естественное преобразование $\Delta (f):\Delta (X)\to \Delta (Y)$, сопоставляющее каждой вершине графа $G$ стрелку $f:X\to Y$.
\end{itemize}
\end{proposition}
Доказательство предоставляется читателю.

\begin{definition}
Диаграмма $D:G\to {\cal K}$ имеет {\bf предел} $A\in {\cal K}_0$, если контравариантный функтор $\mathbf {Dgrm}_{G,{\cal K}}(\Delta (-),D):=\mathbf {Dgrm}_{G,{\cal K}}(-,D)\circ \Delta :{\cal K}^{op}\to \mathbf {Set}$ является представимым с представляющим объектом $A$.  
\end{definition}

\noindent Другими словами, диаграмма $D:G\to {\cal K}$ имеет предел $A\in {\cal K}_0$, если существует естественная биекция $\rho :{\cal K}(-,A)\stackrel {\simeq }{\to }\mathbf {Dgrm}_{G,{\cal K}}(\Delta (-),D)$. Универсальный элемент $\rho _A(1_A)\in \mathbf {Dgrm}_{G,{\cal K}}(\Delta (A),D)$ называется {\bf универсальным конусом} над диаграммой $D$ с {\bf вершиной} $A$.

\addvspace{0.2cm}
\noindent Фиксируем какой-нибудь объект $B\in {\cal K}_0$. Естественные преобразования в $\mathbf {Dgrm}_{G,{\cal K}}(\Delta (B),D)$ состоят из стрелок с началом $B$ и концом в $D$, по одной для каждого объекта в $D$ так, что все получающиеся треугольники с вершиной $B$ и одной из сторон, лежащей в $D$, коммутируют. По формуле $\mu =F(-)(a)$ получаем выражение для $\rho $. Запишем $B$-компоненту естественного преобразования $\rho $. $\rho _B=(\mathbf {Dgrm}_{G,{\cal K}}(\Delta (-),D)(\rho _A(1_A)))_B:{\cal K}(B,A)\to \mathbf {Dgrm}_{G,{\cal K}}(\Delta (B),D):f\mapsto \rho _A(1_A)\circ \Delta (f)$, что сводится к простому умножению справа на $f$ стрелок универсального конуса так, что получается новый конус с вершиной $B$. Рассмотрим для примера рисунок\hskip 0.2cm
$\vcenter
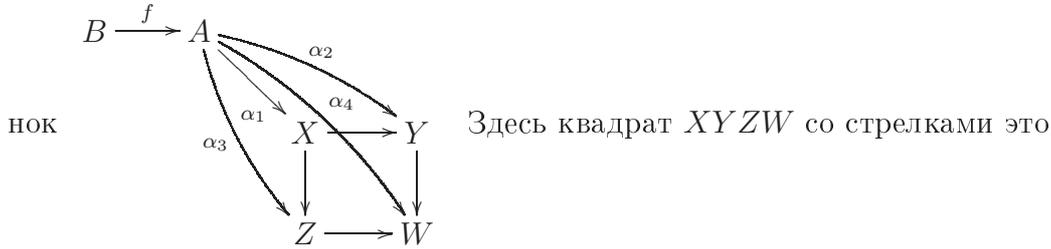
{\xymatrix{
B \ar[r]^-{f} & A \ar[dr]_(0.7){\alpha _1} \ar@/^0.5pc/[drr]^-{\alpha _2} \ar@/_0.5pc/[ddr]_-{\alpha _3} \ar@/^0.65pc/[ddrr]^-{\alpha _4} &   & \\
  &   & X \ar[r] \ar[d] & Y \ar[d] \\
  &   & Z \ar[r] & W
}}$ 
\hskip 0.2cmЗдесь квадрат $XYZW$ со стрелками это диаграмма $D$, $\{\alpha _{i}\},\ i=1,..,4$, это универсальный конус над $D$, $f$ определяет естественное преобразование между постоянными диаграммами $B$ и $A$, $\{\alpha _i\circ f\},\ i=1,..,4$ -- новый конус с вершиной $B$.

\addvspace{0.2cm}
\noindent По условию $\rho _B$ является биекцией. Это значит, что для любого конуса $\beta $ над диаграммой $D$ с вершиной $B$ (или, что то же самое, для любого естественного преобразования $\beta :\Delta (B)\to D$) существует {\bf единственная} стрелка $f:B\to A$ такая, что $\beta =\rho _A(1_A)\circ \Delta (f)$ (то есть, говорят, что конус $\beta $ единственным образом {\bf пропускается} через универсальный конус $\alpha :=\rho _A(1_A)$). Другими словами, универсальный конус является {\bf конечным} объектом в категории конусов над диаграммой $D$, и по этой причине определен однозначно с точностью до изоморфизма. Формулировка того, что такое категория конусов над диаграммой, остается читателю.

\addvspace{0.1cm}
Ниоткуда не следует, что предел той или иной диаграммы существует. Действительно, не все диаграммы в произвольной категории имеют предел. Категория, в которой все диаграммы имеют предел, называется {\bf полной} (в смысле пределов), а функтор, сохраняющий пределы, {\bf непрерывным}. 

\addvspace{0.2cm}
{\bf Примеры.}
\begin{enumerate}
\item Конечный объект  $\mathbf 1$ категории ${\cal K}$ рассматривают как предел пустой диаграммы. Пустая диаграмма это единственный объект категории $\mathbf {Dgrm}_{\O ,{\cal K}}$. В этом случае все (в том числе постоянные) диаграммы совпадают с пустой диаграммой, а все естественные преобразования между диаграммами совпадают с тождественным пустым преобразованием пустой диаграммы. Поэтому, ${\cal K}(X,\mathbf 1)\stackrel {\simeq }{\to }\mathbf {Dgrm}_{\O ,{\cal K}}(\Delta (X),D)$ (изоморфизм между одноэлементными множествами, естественный в $X$).
\item Предел диаграммы, состоящей из двух объектов $X \hskip 0.5cmY$ (для графа $G=\bullet \hskip 0.5cm\bullet $), называется {\bf произведением} (или иногда прямым или декартовым произведением) объектов $X$ и $Y$, и обозначается $X\times Y$. Универсальный конус состоит из двух стрелок $\vcenter{\xymatrix{
 & X\times Y \ar[dl]_-{pr_1} \ar[dr]^-{pr_2} & \\
X &  & Y
}}$ 
которые называются первой и второй {\bf проекциями}. Произведение характеризуется тем, что для любых двух стрелок $g:Z\to X$ и $h:Z\to Y$ существует единственная стрелка $f:Z\to X\times Y$ такая, что диаграмма
$\vcenter{\xymatrix{
 & Z \ar@{-->}[d]_-{\exists !\, f} \ar[dl]_-{g} \ar[dr]^-{h} & \\
   X & X\times Y \ar[l]^-{pr_1} \ar[r]_-{pr_2} & Y 
}}$  
коммутативна. 

\noindent Если в какой-то категории произведение окажется построенным двумя различными способами (что иногда случается), то между такими объектами существует {\it единственный} изоморфизм, переводящий один универсальный конус в другой.

\addvspace{0.1cm}
\noindent В категории $\mathbf {Set}$ произведение двух множеств определяется стандартно как $X\times Y:=\{(x,y)\, |\, x\in X,\ y\in Y\}$, а проекции как $pr_1(x,y):=x$, $pr_2(x,y):=y$.

\addvspace{0.1cm}
\noindent В категории $\mathbf {Cat}$ произведение двух категорий ${\cal K}$ и ${\cal L}$ является категорией ${\cal K}\times {\cal L}$ с объектами -- парами объектов $(K,L),\ K\in {\cal K}_0,\ L\in {\cal L}_0$ и стрелками -- парами стрелок $(f,g),\ f\in {\cal K}_1,\ g\in {\cal L}_1$ с покомпонентным умножением $(f_2,g_2)\circ (f_1,g_1):=(f_2\circ f_1,g_2\circ g_1)$. Проекциями являются функторы $pr_1:{\cal K}\times {\cal L}\to {\cal K}$ и $pr_2:{\cal K}\times {\cal L}\to {\cal L}$, выделяющие соответственно первую и вторую компоненты объектов или стрелок. 

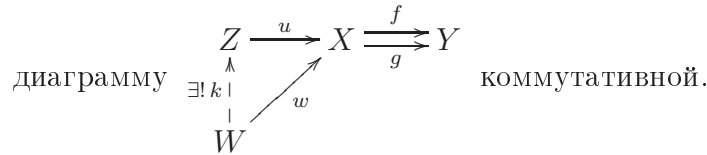
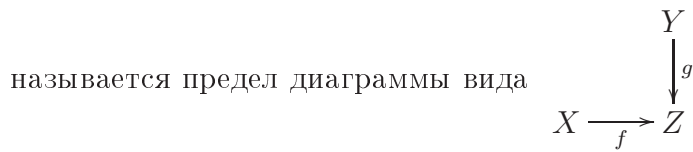
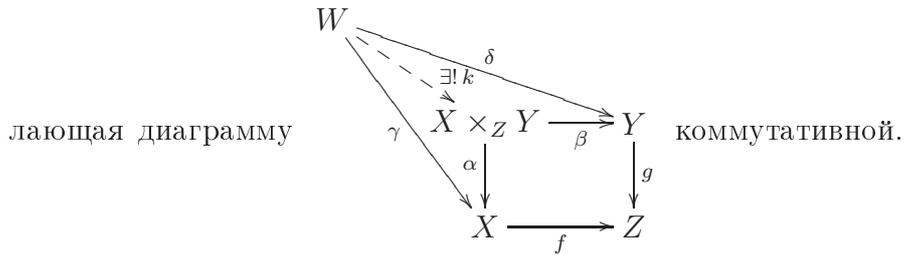
\addvspace{0.1cm}
Произведение трех или большего числа объектов определяется аналогично как предел диаграммы, состоящей из $n$ объектов. Читателю рекомендуется доказать, например, для трех объектов, что $3$-кратное произведение объектов $X\times Y \times Z$ и объекты $X\times (Y\times Z)$ или $(X\times Y)\times Z$, полученные итерацией бинарного произведения, изоморфны.   
\item {\bf Уравнитель} это предел диаграммы, состоящей из двух параллельных стрелок $\xymatrix{X \ar@<0.5ex>[r]^-{f} \ar@<-0.5ex>[r]_-{g} & Y}$. Это значит, существует предельный объект $Z$ и универсальный конус $\vcenter{\xymatrix{
Z \ar[dr]_-{u} \ar@/^0.5pc/[drr]^-{v} & & \\
  &  X \ar@<0.5ex>[r]^(0.45){f} \ar@<-0.5ex>[r]_(0.45){g} & Y
}}$
такие, что $v=f\circ u=g\circ u$. 

\addvspace{0.1cm}
\noindent Поскольку компонента $v$ для универсального конуса (и для всякого другого конуса) полностью определена компонентой $u$, то ее (лишнюю компоненту) обычно не рассматривают, и называют уравнителем стрелку $u:Z\to X$ такую, что 
\begin{itemize}
\item $f\circ u=g\circ u$, 
\item для всякой другой стрелки $w:W\to X$, такой, что $f\circ w=g\circ w$, существует единственная стрелка $k:W\to Z$, делающая диаграмму 
$\vcenter{\xymatrix{
Z \ar[r]^-{u} & X \ar@<0.5ex>[r]^-{f} \ar@<-0.5ex>[r]_-{g} & Y\\
W \ar[ur]_-{w} \ar@{-->}[u]^-{\exists !\, k} & & 
}}$ 
коммутативной. 

\addvspace{0.1cm}
\noindent  Как и все пределы, уравнитель, если существует, определен однозначно с точностью до изоморфизма.

\addvspace{0.1cm}
В категории множеств $\mathbf {Set}$ уравнителем $u$ будет вложение подмножества $\{x\in X\, |\, f(x)=g(x)\}$, на котором $f$ и $g$ совпадают, в множество $X$. Аналогично, в категории $\mathbf {Vect}$ уравнителем будет вложение подпространства $\{x\in X\, |\, f(x)=g(x)\}$, на котором линейные отображения $f$ и $g$ совпадают, в векторное пространство $X$. 
\end{itemize}
\item {\bf Расслоенным произведением} $X\times _{Z}Y$ объектов $X$ и $Y$ над $Z$ называется предел диаграммы вида\hskip 0.2cm $\vcenter{\xymatrix{
   &  Y \ar[d]^-{g}\\
X \ar[r]_-{f} & Z
}}$ 

\addvspace{0.2cm}
Как и в предыдущем примере, стрелка конуса в объект $Z$ вычисляется через остальные стрелки конуса над указанной диаграммой и не учитывается. Таким образом, расслоенным произведением называется объект $X\times _{Z}Y$ и пара стрелок $\alpha :X\times _{Z}Y\to X$, $\beta :X\times _{Z}Y\to Y$ таких, что
\begin{itemize}
\item $f\circ \alpha =g\circ \beta $,
\item для любых стрелок $\gamma :W\to X$, $\delta :W\to Y$ с условием $f\circ \gamma =g\circ \delta $ существует единственная стрелка $k:W\to X\times _{Z}Y$, делающая диаграмму $\vcenter{\xymatrix{
W \ar[ddr]_-{\gamma } \ar[drr]^-{\delta } \ar@{-->}[dr]^(0.7){\exists !\, k} & &  \\
  & X\times _{Z}Y \ar[d]_(0.4){\alpha } \ar[r]_-{\beta } & Y \ar[d]^-{g} \\
  & X \ar[r]_-{f}  &  Z
}}$
коммутативной.
Квадрат со сторонами $\alpha ,f,\beta ,g$ называется {\bf универсальным квадратом}.
\end{itemize} 
В категории множеств $\mathbf {Set}$ расслоенным произведением $X\times _{Z}Y$ будет подмножество декартова произведения $X\times Y$, на котором функции $f\circ pr_1$ и $g\circ pr_2$ совпадают, то есть $X\times _{Z}Y=\{(x,y)\in X\times Y\, |\, f(x)=g(y)\}$. Функции $\alpha $ и $\beta $ являются соответственно ограничениями проекций $pr_1:X\times Y\to X$ и $pr_2:X\times Y\to Y$ на данное подмножество. 
\end{enumerate}

Приведенные примеры представляют собой в некотором смысле элементарные кирпичики, из которых строятся более сложные пределы. Имеется теорема, что категория является конечно-полной (то есть пределы всех диаграмм, состоящих из конечного числа объектов и стрелок, существуют), если существуют конечный объект, бинарные произведения и уравнители (или если существуют конечный объект и расслоенные произведения). 

\begin{assignment}
\begin{enumerate}
\item Доказать утверждение 9.
\item Нарисовать коммутативные диаграммы, которые должны выполняться, для стрелок конуса над диаграммой.
\item Для категории ${\cal K}$ и объекта $K\in {\cal K}_0$ составим новую категорию, которая обозначается ${\cal K}/K$ и  называется категорией стрелок над объектом $K$. Ее объектами являются все стрелки с концом $K$, а морфизмами из объекта $\vcenter{\xymatrix{X \ar[d]_-{f}\\ K}}$ в объект $\vcenter{\xymatrix{Y \ar[d]_-{g}\\ K}}$ стрелки категории ${\cal K}$ с началом $X$ и концом $Y$, для которых треугольник $\vcenter{\xymatrix{X \ar[rr] \ar[dr]_-{f} && Y \ar[dl]^-{g} \\
  & K & }}$ коммутативен. Доказать, что эта конструкция действительно определяет категорию. Доказать, что категория конусов над диаграммой $D$ представляет собой полную подкатегорию категории стрелок $\mathbf {Dgrm}_{G,{\cal K}}/D$. Убедиться, что универсальный конус над $D$ является конечным объектом в этой подкатегории.
\item Разобрать случай предела пустой диаграммы (пример 1).
\item Проверить, что произведения в категориях множеств $\mathbf {Set}$ и категорий $\mathbf {Cat}$ определены корректно (пример 2).
\item Найти произведение в категории векторных пространств $\mathbf {Vect}$ и предупорядоченном множестве $(L,<)$.
\item Доказать, что $X\times Y\times Z\simeq X\times (Y\times Z)\simeq (X\times Y)\times Z$.
\item Доказать, что уравнитель является монострелкой.
\item Определить смысл универсального квадрата в категории $\mathbf {Set}$, для которого $f$ -- произвольная функция, а $g$ -- вложение подмножества $Y$ в множество $Z$ (пример 4).
\item В категории ${\cal K}$ с бинарными произведениями определить функтор умножения $X\times (-):{\cal K}\to {\cal K}$ на фиксированный объект категории $X\in {\cal K}_0$. 
\end{enumerate}
\end{assignment}

\section{Копределы}

Копределом диаграммы является предел двойственной диаграммы в двойственной категории. Двойственная диаграмма, как и двойственная категория, получается заменой всех стрелок на противоположные. Имеется следующая последовательность изоморфизмов ${\cal K}(A,-)\simeq {\cal K}^{op}(-,A)\simeq \mathbf {Dgrm}_{G,{\cal K}^{op}}(\Delta (-),D^{op})\simeq \mathbf {Dgrm}_{G,{\cal K}}(D,\Delta (-))$.

\begin{definition}
Объект $A\in {\cal K}_0$ называется {\bf копределом} диаграммы $D:G\to {\cal K}$, если существует естественный изоморфизм $\rho :{\cal K}(A,-)\stackrel {\simeq }{\to }\mathbf {Dgrm}_{G,{\cal K}}(D,\Delta (-))$. Универсальный элемент $\rho (1_A)$ называется {\bf универсальным коконусом} над диаграммой $D$.
\end{definition}

Графически коконус над диаграммой выглядит как семейство стрелок, начинающихся в вершинах диаграммы $D$ и заканчивающихся в одном объекте -- {\bf вершине коконуса}. {\it Универсальный коконус} является начальным объектом в категории коконусов над диаграммой, то есть любой коконус пропускается через него единственным образом.

\addvspace{0.2cm}
{\bf Примеры.}
\begin{enumerate}
\item Начальный объект $\mathbf 0$ является копределом пустой диаграммы $D:\O \to {\cal K}$. 
\item {\bf Копроизведением} $X\coprod Y$ называется копредел диаграммы, состоящей из двух объектов $X$ и $Y$. Универсальный коконус состоит из двух стрелок $\vcenter{\xymatrix{X \ar[r]^-{i_1} & X\coprod Y & Y \ar[l]_-{i_2} }}$, называемых {\bf вложениями}. Любой коконус над диаграммой \, $X \hskip 0.5cm Y$ \, единственным образом пропускается через универсальный коконус \hskip -0.3cm$\vcenter{\xymatrix{
   &   Z   &     \\
X \ar[r]_-{i_1} \ar[ur]^-{f} & X\coprod Y \ar@{-->}[u]^-{\exists !\, k} & Y  \ar[l]^-{i_2} \ar[ul]_-{g}
}}$

\noindent В категории множеств $\mathbf {Set}$ копроизведением является дизъюнктное объединение множеств $X\coprod Y:=(X\times \{0\})\bigcup (Y\times \{1\})$ (элементы множества $X$ берутся с меткой $0$, а множества $Y$ с меткой $1$, чтобы быть уверенным, что множества не пересекаются. Дизъюнктное объединение это объединение непересекающихся множеств).
\item {\bf Коуравнителем} называется копредел диаграммы, состоящей из двух параллельных стрелок 
$\xymatrix{X \ar@<0.5ex>[r]^-{f} \ar@<-0.5ex>[r]_-{g} & Y}$. Как и в случае пределов, коконус определяется одной стрелкой, выходящей из объекта $Y$. Таким образом, коуравнитель это объект $Z$ вместе со стрелкой $u:Y\to Z$ такими, что
\begin{itemize}
\item $u\circ f=u\circ g$,
\item для любой стрелки $w:Y\to W$ со свойством $w\circ f=w\circ g$ существует единственная стрелка $k:Z\to W$, делающая диаграмму
$\vcenter{\xymatrix{
 &   &  W \\
X \ar@<0.5ex>[r]^-{f} \ar@<-0.5ex>[r]_-{g} & Y \ar[r]_-{u} \ar[ur]^-{w}  & Z \ar@{-->}[u]_-{\exists !\, k}
}}$
коммутативной.
\end{itemize}

\addvspace{0.2cm}
\noindent В категории множеств $\mathbf {Set}$ коуравнители определяются посложнее, чем уравнители. Именно, $Z$ определяется как множество классов эквивалентности по отношению эквивалентности, {\it порожденному} условиями $y_1\sim y_2$, где $y_1,y_2\in Y$, если существует элемент $x\in X$ такой, что $f(x)=y_1$, $g(x)=y_2$. Другими словами, значения функций $f(x)$ и $g(x)$ считаются эквивалентными в каждой точке $x\in X$, и строится (минимальное) отношение эквивалентности по этим условиям так, чтобы, например, выполнялось $y\sim y$, $f(x)\sim f(x)$, $g(x)\sim f(x)$, $(f(x)\sim g(x))\& (f(x')\sim g(x'))\& (f(x)=f(x'))\Rightarrow (g(x)\sim g(x'))$, и так далее. Если не строить отношение эквивалентности, то нельзя получить множество $Z$. Универсальной стрелкой $u:Y\to Z$ в этом случае будет каноническая проекция, сопоставляющая элементу $y\in Y$ его класс эквивалентности. 
\item {\bf Приклеиванием} объекта $Y$ к объекту $X$ по объекту $Z$ называется копредел диаграммы, состоящей из двух стрелок с общим началом $\vcenter{\xymatrix{Z \ar[r]^-{f} \ar[d]_-{g} & X \\
Y & }}$
Этот копредел часто обозначается $X\coprod _{Z}Y$ и может называться дизъюнктным объединением объектов $X$ и $Y$ над объектом $Z$. Коконусы над данной диаграммой определены двумя стрелками выходящими из вершин $X$ и $Y$. Таким образом, {\it дизъюнктное объединение объектов $X$ и $Y$ над объектом $Z$} состоит из объекта $X\coprod _{Z}Y$ вместе с двумя стрелками 
$\vcenter{\xymatrix{
  & X \ar[d]^-{\alpha } \\
Y \ar[r]_-{\beta } & X\coprod _{Z}Y
}}$  
такими, что 
\begin{itemize}
\item $\alpha \circ f=\beta \circ g$,
\item для любых двух стрелок $\gamma :X\to W$, $\delta :Y\to W$ со свойством $\gamma \circ f=\delta \circ g$ существует единственная стрелка $k:X\coprod _{Z}Y\to W$, делающая диаграмму 
$\vcenter{\xymatrix{
Z \ar[r]^-{f} \ar[d]_-{g} & X \ar[d]_-{\alpha } \ar[ddr]^-{\gamma } &  \\
Y \ar[r]^-{\beta } \ar[drr]_-{\delta } & X\coprod _{Z}Y \ar@{-->}[dr]^(0.3){\exists !\, k} \\
    &       &     W
}}$
коммутативной. 

\addvspace{0.2cm}
Квадрат $f, \alpha , g, \beta $ называется {\bf коуниверсальным квадратом}.
\end{itemize}

\addvspace{0.2cm}
\noindent В категории множеств $\mathbf {Set}$ результат приклеивания можно описать как множество классов эквивалентности дизъюнктного объединения множеств $X\coprod Y:=X\times \{0\}\bigcup Y\times \{1\}$ по отношению эквивалентности $\sim $, порожденному парами $i_1\circ f(z)\sim i_2\circ g(z)$, $z\in Z$, где $i_1$ и $i_2$ вложения множеств $X$ и $Y$ в дизъюнктное объединение $X\coprod Y$.
\end{enumerate} 

\addvspace{0.1cm}
Категория, в которой все копределы существуют, называется {\bf кополной}, а если существуют только копределы конечных диаграмм (диаграмм с конечным числом объектов и стрелок), то {\bf конечно кополной}. Функтор, сохраняющий копределы, называется {\bf конепрерывным}. Как и в случае пределов, для того чтобы категория была конечно кополной необходимо и достаточно, чтобы существовали начальный объект, бинарные копроизведения и коуравнители (или же начальный объект и коуниверсальные квадраты). 

\begin{assignment}
\begin{enumerate}
\item Выписать подробно, что означает определение копредела диаграммы, как это было сделано в параграфе 2.2 для предела.
\item Доказать, что дизъюнктное объединение множеств является копроизведением в категории множеств $\mathbf {Set}$.
\item Убедиться в том, что любое бинарное отношение $\rho $ на множестве $A$ порождает некоторое отношение эквивалентности $\sim$ на $A$, содержащее $\rho $ (Бинарное отношение $\rho $ это, по определению, любое подмножество $A^2:=A\times A$. Всегда существуют отношения эквивалентности (то есть рефлексивные, симметричные и транзитивные отношения), содержащие отношение $\rho $, например, $\sim \, :=A^2$. Пересечение любого числа отношений эквивалентности является снова отношением эквивалентности).
\item Доказать, что коуравнитель в категории $\mathbf {Set}$ корректно определен (пример 3).
\item Убедиться, что если $Z\subset X$ и $Z\subset Y$ в категории множеств $\mathbf {Set}$, то $X\coprod _{Z}Y$ действительно имеет смысл 'приклеивания' множества $Y$ к множеству $X$ по точкам подмножества $Z$. Рекомендуется также проверить это в категории топологических пространств $\mathbf {Top}$.
\end{enumerate}
\end{assignment}

\section{Сопряженные функторы}

Сопряженные функторы встречаются в разных областях математики и имеют важное значение. Часто какой-нибудь тривиальный функтор (например, забывающий или включения подкатегории) имеет нетривиальный сопряженный. На этом пути устанавливаются интересные взаимоотношения между категориями, некоторые понятия (например, {\it свободной} алгебраической системы) находят свое обоснование, вводятся дополнительные структуры в категории, и так далее. С некоторой точки зрения аппарат сопряженных функторов является основным в теории категорий.

\begin{definition}
Пусть функторы $F:{\cal L}\to {\cal K}$ и $H:{\cal K}\to {\cal L}$ действуют в противоположных направлениях. Тогда $F$ называется {\bf левым сопряженным} (функтору $H$), а $H$ {\bf правым сопряженным} (функтору $F$), если существует естественная биекция $\theta _{X,Y}:{\cal K}(F(X),Y)\stackrel {\simeq }{\to }{\cal L}(X,H(Y))$, где $X\in {\cal L}_0$, $Y\in {\cal K}_0$.
\end{definition} 

\noindent Ситуация сопряжения обозначается $F\dashv H$ ($F$ сопряжен слева $H$).

\addvspace{0.2cm}
\noindent В определении говорится, что существует взаимно-однозначное соответствие между {\it стрелками}, принадлежащими $Hom$-множествам двух разных категорий ${\cal K}$ и ${\cal L}$. Естественность в аргументах $X$ и $Y$ означает, что для любой ${\cal L}$-стрелки $f:X'\to X$ и для любой ${\cal K}$-стрелки $g:Y\to Y'$ диаграмма
$\vcenter{\xymatrix{
X' \ar[d]_-{f} &  Y'              &  {\cal K}(F(X'),Y') \ar[r]^-{\theta _{X',Y'}}  &   {\cal L}(X',H(Y')) \\
X              &  Y  \ar[u]^-{g}  &  {\cal K}(F(X),Y) \ar[r]_-{\theta _{X,Y}}  \ar[u]^-{{\cal K}(F(f),g)} &   {\cal L}(X,H(Y)) \ar[u]_-{{\cal L}(f,H(g))}
}}$
коммутативна, то есть для любых объектов $X\in {\cal L}_0$ и $Y\in {\cal K}_0$, и для любой стрелки $x:F(X)\to Y$ выполняется тождество $\theta _{X',Y'}(g\circ x\circ F(f))=H(g)\circ \theta _{X,Y}(x)\circ f$.

\addvspace{0.2cm}
\noindent Стрелки $x:F(X)\to Y$ и $\theta _{X,Y}(x):X\to H(Y)$ называются {\bf сопряженными} друг другу, и переход от одной к другой часто обозначается чертой сверху ($\overline{x}:=\theta _{X,Y}(x)$, $x:=\overline{\theta _{X,Y}(x)}$, так что $\overline{\overline{x}}=x$).  Графически сопряженные стрелки изображают также ${\displaystyle \frac{F(X)\stackrel {x}{\to }Y}{X\stackrel {\overline{x}}{\to }H(Y)}}$. Тогда естественность операции $\overline{(\cdot )}$ (или, что то же самое, естественность $\theta $) будет выглядеть как ${\displaystyle \frac{F(X')\stackrel {F(f)}{\to }F(X)\stackrel {x}{\to }Y\stackrel {g}{\to }Y'}{X'\stackrel {f}{\to }X\stackrel {\overline{x}}{\to }H(Y)\stackrel {H(g)}{\to }H(Y')}}$.

\addvspace{0.2cm}
\noindent Из формулы $\theta _{X,Y}:{\cal K}(F(X),Y)\stackrel {\simeq }{\to } {\cal L}(X,H(Y))$ следует, что все функторы ${\cal K}(F(-),Y):{\cal L}^{op}\to \mathbf {Set}$ (зависящие от параметра $Y\in {\cal K}_0$) представимы с представляющим объектом $H(Y)\in {\cal L}_0$ и универсальным элементом $\theta ^{-1}_{HY,Y}(1_{HY})$. Точно также все функторы ${\cal L}(X,H(-)):{\cal K}\to \mathbf {Set}$ (параметризованные объектами $X\in {\cal L}_0$) представимы с представляющим объектом $F(X)\in {\cal K}_0$ и универсальным элементом $\theta _{X,FX}(1_{FX})$.

\addvspace{0.2cm}
\noindent Универсальный элемент $\theta ^{-1}_{HY,Y}(1_{HY}):FHY\to Y$ называется {\bf коединицей} (ситуации) сопряжения (на объекте $Y$) и обозначается $\varepsilon _{Y}:FHY\to Y$, а универсальный элемент $\theta _{X,FX}(1_{FX}):X\to HFX$ соответственно называется {\bf единицей} (ситуации) сопряжения (на объекте $X$) и обозначается $\eta _{X}:X\to HFX$. Таким образом, ${\displaystyle \frac{FHY\stackrel {\varepsilon _Y}{\to }Y}{HY\stackrel {1_{HY}}{\to }HY}}$ и ${\displaystyle \frac{FX\stackrel {1_{FX}}{\to }FX}{X\stackrel {\eta _X}{\to }HFX}}$.

\begin{proposition}
\begin{itemize}
\item Стрелки $\varepsilon _Y:FHY\to Y$ являются компонентами естественного преобразования $\varepsilon :FH\Rightarrow 1_{\cal K}$.
\item Стрелки $\eta _X:X\to HFX$ являются компонентами естественного преобразования $\eta :1_{\cal L}\Rightarrow HF$. 
\end{itemize}
\end{proposition}

Доказательство.
Докажем только первое утверждение.

\noindent Рассмотрим диаграмму $\vcenter{\xymatrix{
FHY \ar[r]^-{\varepsilon _Y} & Y \\
FHY' \ar[r]_-{\varepsilon _{Y'}} \ar[u]^-{FHf} & Y' \ar[u]_-{f}
}}$ которая предполагается быть коммутативной. \hskip 0.2cmВозьмем сопряженные стрелки обоих произведений

\addvspace{0.2cm}
\noindent ${\displaystyle \frac{FHY'\stackrel {FHf}{\to }FHY\stackrel {\varepsilon _Y}{\to }Y}{HY'\stackrel {Hf}{\to }HY\stackrel {1_{HY}}{\to }HY}}$ \hskip 0.2cmи \hskip 0.2cm 
${\displaystyle \frac{FHY'\stackrel {\varepsilon _{Y'}}{\to }Y'\stackrel {f}{\to }Y}{HY'\stackrel {1_{HY'}}{\to }HY'\stackrel {Hf}{\to }HY}}$. \hskip 0.2cmОни совпадают. 

\addvspace{0.2cm}
\noindent Следовательно, исходные стрелки тоже совпадают.   \hfill  $\square $

\begin{proposition}
Следующие утверждения эквивалентны
\begin{itemize}
\item $F\dashv H$,
\item для каждого $Y\in {\cal K}_0$ функтор ${\cal K}(F(-),Y):{\cal L}^{op}\to \mathbf {Set}$ представим, 
\item для каждого $X\in {\cal L}_0$ функтор ${\cal L}(X,H(-)):{\cal K}\to \mathbf {Set}$ представим.
\end{itemize}
\end{proposition}

Доказательство.
То, что из первого утверждения следуют остальные, было рассмотрено в начале параграфа. Докажем, что из второго пункта следует первый. Представимость функтора ${\cal K}(F(-),Y):{\cal L}^{op}\to \mathbf {Set}$ означает, что для всякого объекта $Y\in {\cal K}_0$ существует объект $HY\in {\cal L}_0$ такой, что ${\cal K}(F(-),Y)\simeq {\cal L}(-,HY)$. Требуется доказать: что отображение $H$ может быть распространено функториально над стрелками в ${\cal K}$, определить биекцию $\theta _{X,Y}$ и доказать ее естественность.

\noindent Прежде всего, из представимости функтора ${\cal K}(F(-),Y):{\cal L}^{op}\to \mathbf {Set}$ с представляющим объектом $HY$ и универсальным элементом $\varepsilon _{Y}:FHY\to Y$ следует, что 'категория элементов' (смотрите параграф 2.1), построенная по функтору ${\cal K}(F(-),Y)$, имеет пару $(HY,FHY\stackrel {\varepsilon _{Y}}{\to }Y)$ своим конечным объектом, то есть для каждой другой пары $(X,FX\stackrel {x}{\to }Y)$ существует единственная стрелка $\overline {x}:X\to HY$ такая, что диаграмма $\vcenter{\xymatrix{
HY  & FHY  \ar[r]^-{\varepsilon _{Y}} & Y \\
X \ar@{-->}[u]^-{\exists !\, \overline x}  &  FX \ar[ur]_-{\forall x} \ar[u]^-{F(\overline {x})} & 
}}$
коммутативна. Отсюда следуют все требуемые 

\vspace{0.3cm}
\noindent результаты.

\addvspace{0.2cm}
\noindent Определим функтор $H$ на стрелках следующим образом $Hf:=\overline {f\circ \varepsilon _{Y'}}$, то есть
$\vcenter{\xymatrix{
HY              &     FHY  \ar[r]^-{\varepsilon _Y}  & Y  \\
HY' \ar@{-->}[u]^-{\exists !\, Hf} & FHY' \ar[u]^-{FHf} \ar[r]_-{\varepsilon _{Y'}} & Y' \ar[u]_-{f}
}}$\hskip 0.3cm 
(то, что $H$ функтор, то есть сохраняет единицы и произведения, остается доказать читателю)

\addvspace{0.2cm}
\noindent   
Из последней диаграммы следует также, что $\varepsilon :FH\Rightarrow 1_{\cal K}$ есть естественное преобразование.

\addvspace{0.2cm}
\noindent Определим $\theta _{X,Y}:{\cal K}(FX,Y)\to {\cal L}(X,HY):x\mapsto \overline x$. Очевидно, это биекция с обратным отображением $\theta ^{-1}_{X,Y}:{\cal L}(X,HY)\to {\cal K}(FX,Y):y\mapsto \varepsilon _{Y}\circ F(y)$.

\addvspace{0.2cm}
\noindent Рассмотрим диаграмму $\vcenter{\xymatrix{
X' \ar[d]_-{f} & Y' & {\cal K}(FX',Y') \ar[r]^-{\theta _{X',Y'}} & {\cal L}(X',HY') \\
X  & Y \ar[u]^-{g} &  {\cal K}(FX,Y) \ar[r]_-{\theta _{X,Y}} \ar[u]^-{{\cal K}(Ff,g)} & {\cal L}(X,HY) \ar[u]_-{{\cal L}(f,Hg)}
}}$

\addvspace{0.1cm}
\noindent Она коммутативна, потому что сопряженные двух возможных произведений совпадают. Действительно, для всякой стрелки $x:FX\to Y$ имеем $\theta ^{-1}_{X',Y'}(\theta _{X',Y'}(g\circ x\circ Ff))=g\circ x\circ Ff$ \, и \, $\theta ^{-1}_{X',Y'}(Hg\circ \theta _{X,Y}(x) \circ f)=\varepsilon _{Y'}\circ F(Hg\circ \overline x\circ f)=\varepsilon _{Y'}\circ  FHg\circ F\overline x\circ Ff=g\circ \varepsilon _Y\circ F\overline x\circ Ff=g\circ x\circ Ff$.     \hfill $\square $ 

\addvspace{0.2cm}
{\bf Примеры.}
\begin{enumerate}
\item Вспомним определения предела и копредела:

\noindent ${\cal K}(-,A)\stackrel {\simeq }{\to }\mathbf {Dgrm}_{G,{\cal K}}(\Delta (-),D)$ \, и 

\noindent ${\cal K}(A',-)\stackrel {\simeq }{\to }\mathbf {Dgrm}_{G,{\cal K}}(D,\Delta (-))$. 

Обозначим операции, назначающие диаграмме $D:G\to {\cal K}$ предельный и копределеный объекты, через $\text{\tt lim}$ и $\text{\tt colim}$. То есть $A=\text{\tt lim}\, D$, $A'=\text{\tt colim}\, D$ (предполагается, что пределы и копределы существуют в ${\cal K}$). Функции $\text{\tt lim}$ и $\text{\tt colim}$ однозначно продолжаются до функторов $\text{\tt lim},\text{\tt colim}:\mathbf {Dgrm}_{G,{\cal K}}\to {\cal K}$ (доказательство этого факта выносится в упражнения). Следовательно, имеем 

\noindent ${\cal K}(-,\text{\tt lim}\, D)\stackrel {\simeq }{\to }\mathbf {Dgrm}_{G,{\cal K}}(\Delta (-),D)$ \, и

\noindent ${\cal K}(\text{\tt colim}\, D,-)\stackrel {\simeq }{\to }\mathbf {Dgrm}_{G,{\cal K}}(D,\Delta (-))$, \, то есть

\noindent $\text{\tt colim}\dashv \Delta \dashv \text{\tt lim}$ (копредел является левым сопряженным, а предел правым сопряженным постоянному функтору).
\item Забывающий функтор $U:\mathbf {Mon}\to \mathbf {Set}$ является правым сопряженным к функтору взятия {\bf свободного моноида} $F:\mathbf {Set}\to \mathbf {Mon}$, который сопоставляет множеству $A$ моноид 'слов', для которых $A$ является 'алфавитом'. Словами считаются любые конечные последовательности элементов $A$, например $a_1a_2a_3$ или $abcd$ (если все буквы принадлежат множеству $A$). Единицей считается пустое слово, а умножением операция приписывания одного слова к другому, например, $(a_1a_2a_3)\cdot (abcd)=a_1a_2a_3abcd$. Отображению $f:A\to B$ функтор $F$ сопоставляет гомоморфизм моноидов $Ff:FA\to FB:a_1\dots a_n\mapsto f(a_1)\dots f(a_n)$. Ситуация сопряжения $\mathbf {Mon}(FA,M)\simeq \mathbf {Set}(A,UM):f\mapsto f|_A$ говорит о том, что гомоморфизмы {\it из} свободного моноида $FA$ в другой моноид $M$ находятся во взаимнооднозначном соответствии с отображениями {\it множества $A$, порождающего свободный моноид} в моноид $M$. Более точно ситуация выражается с использованием единицы $\eta _{A}:A\to UFA:a\mapsto a$ (каждый элемент $a$ вкладывается как слово, состоящее из одной буквы). Именно, универсальное свойство единицы дает следующую диаграмму \hskip 0.6cm $\vcenter{\xymatrix{
A \ar[r]^-{\eta _A} \ar[dr]_-{\forall f}  &   UFA \ar[d]^-{U\overline f} & FA \ar@{-->}[d]^-{\exists !\, \overline f} \\
 & UM & M
}}$ 

\noindent которая называется обычно {\bf универсальным свойством свободного моноида}. Гомоморфизм $\overline f:FA\to M$ определяется как $\overline f(a_1\dots a_n):=f(a_1)\dots f(a_n)$. Легко видеть, что $\overline f$ действительно гомоморфизм и нет другого пути его определения.

\noindent Аналогично, забывающие функторы $\mathbf {Grp}\to \mathbf {Set}$, $\mathbf {Ab}\to \mathbf {Set}$, $\mathbf {Vect}\to \mathbf {Set}$ и другие имеют левые сопряженные (функторы взятия свободного объекта).   
\item Забывающий функтор $U:\mathbf {Top}\to \mathbf {Set}$ является одновременно правым и левым сопряженным. Левый сопряженный к нему $F:\mathbf {Set}\to \mathbf {Top}$ сопоставляет множеству $A$ топологическое пространство, состоящее из множества $A$, в котором все подмножества считаются открытыми (так называемая {\it дискретная топология}). Имеет место естественная биекция $\mathbf {Top}(FA,X)\stackrel {\simeq }{\to }\mathbf {Set}(A,UX):f\mapsto Uf$, которая выражает факт, что каждое отображение множества $A$ в топологическое пространство $X$ непрерывно в дискретной топологии на $A$ (действительно, прообраз открытого множества в $X$ будет некоторым подмножеством в $A$, которое открыто по определению). Правый сопряженный к $U$ (обозначим его $W:\mathbf {Set}\to \mathbf {Top}$) сопоставляет множеству $A$ топологическое пространство, состоящее из множества $A$, в котором открытыми множествами считаются только подмножества $A$ и $\O $ (так называемая {\it антидискретная топология}). Имеет место естественная биекция $\mathbf {Top}(X,WA)\stackrel {\simeq }{\to }\mathbf {Set}(UX,A):f\mapsto Uf$, выражающая факт, что любая функция из топологического пространства $X$ в антидискретное пространство $WA$ непрерывна (потому что прообраз $A$ есть $X$, прообраз $\O $ есть $\O $, и все множество и пустое множество открыты в любой топологии по определению).
\item Ситуацией сопряжения между двумя предупорядоченными множествами будет $\xymatrix{(L,<)  \ar@/^0.5pc/[r]^-{f}_-{\top } &    (L',<) \ar@/^0.5pc/[l]^-{g} }$, где $f$, $g$ -- монотонные функции, и $g(x)<y$ в $L$, если и только если $x<f(y)$ в $L'$. Единицей сопряжения будет $x<fg(x)$, а коединицей $gf(y)<y$. Сопряжение между предупорядоченными множествами называется {\bf соответствием Галуа}. {\it Обратным соответствием Галуа} называется сопряжение $\xymatrix{(L,<)^{op}  \ar@/^0.5pc/[r]^-{f}_-{\top } &    (L',<) \ar@/^0.5pc/[l]^-{g} }$. В этом случае функции $f$ и $g$ меняют порядок на противоположный (антитонные функции), и $y<f(x)$, если и только если $x<g(y)$.  
\item Ситуация сопряжения не обязана быть между разными категориями. Например, в категории множеств $\mathbf {Set}$ имеется важное сопряжение между функтором $X\times (-)$ умножения на некоторое множество $X$ и функтором $(-)^X:=\mathbf {Set}(X,-)$ взятия множества функций с областью определения $X$. Действительно, имеет место естественный изоморфизм
$\theta _{X,Y}:\mathbf {Set}(X\times Y,Z)\stackrel {\simeq }{\to } \mathbf {Set}(Y,Z^X)$, где $\theta _{X,Y}(f:X\times Y\to Z:(x,y)\mapsto f(x,y)):=\overline f:Y\to Z^X:y\mapsto f(-,y)$. Обратное отображение задается $\theta ^{-1}_{X,Y}(g:Y\to Z^X):=\overline g:X\times Y\to Z:(x,y)\mapsto g(y)(x)$. Проверка естественности $\theta _{X,Y}$ остается читателю. Категории с бинарным произведением, в которых функтор умножения $X\times (-)$ имеет правый сопряженный для каждого объекта $X$, называются {\bf декартово-замкнутыми}. Таким образом, $\mathbf {Set}$ -- декартово-замкнутая категория. $\mathbf {Cat}$ также является декартово-замкнутой, а, например, $\mathbf {Top}$ и $\mathbf {Vect}$ нет.
\end{enumerate}

\addvspace{0.25cm}
Приведем в конце параграфа еще одну теорему без доказательства.

\begin{proposition}
Правые сопряженные являются непрерывными функторами, а левые сопряженные -- конепрерывными.
\end{proposition}

\begin{assignment}
\begin{enumerate}
\item Доказать второй пункт утверждения 10.
\item Доказать, что функция $H$, определенная в доказательстве утверждения 11, является функтором.
\item Доказать, что в утверждении 11 из третьего пункта следует первый.
\item Доказать, что операции взятия предела и копредела однозначно продолжаются до функторов (пример 1).
\item Найти единицу и коединицу сопряжения в примере 1.
\item Определить смысл коединицы в примере 2.
\item Проверить корректность определения левого и правого сопряженных функторов к забывающему функтору $U:\mathbf {Top}\to \mathbf {Set}$ (пример 3).
\item Пусть $A$ есть некоторое множество. Сопоставим каждому его подмножеству $B$ группу преобразований $G_B$ множества $A$, каждый элемент которой оставляет все точки $B$ неподвижными. И наоборот, каждой подгруппе $G$ группы преобразований $\mathbf {Aut}(A)$ множества $A$ сопоставим множество точек в $A$, которые остаются неподвижными при каждом преобразовании подгруппы.  Доказать, что имеет место обратное соответствие Галуа между частично упорядоченными множествами ${\cal P}(A)$ и множеством подгрупп группы $\mathbf {Aut}(A)$. Порядок в каждом множестве задан отношением включения $X\subset Y$. (Группа преобразований $\mathbf {Aut}(A)$ состоит из всех биекций в $\mathbf {Set}(A,A)$).
\item Доказать естественность отображения $\theta _{X,Y}$ в примере 5.
\item Доказать, что произведение левых сопряженных функторов является левым сопряженным функтором, а произведение правых сопряженных -- правым сопряженным функтором. 
\end{enumerate}
\end{assignment}

\section{Сопряженность в 2-категории}

Существует еще одна более абстрактная характеризация сопряженных функторов через тождества, использующие единицу и коединицу. Из этого сразу следует, что 2-функторы на категории $2\text{-}\mathbf {Cat}$ сохраняют сопряжение (потому что все функторы сохраняют категорные тождества). Кроме того это показывает, что сопряженность является чисто алгебраическим понятием и имеет смысл в любой 2-категории.

Рассмотрим две диаграммы \hskip 0.5cm $\vcenter
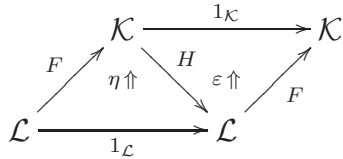
{\xymatrix{
          &  {\cal L} \ar[dr]^(0.44){F} \ar[rr]^-{1_{\cal L}} \ar@{}[d]|{\varepsilon \, \Downarrow } &   \ar@{}[d]|{\eta \, \Downarrow }    &   {\cal L}  \\
{\cal K} \ar[ur]^-{H} \ar[rr]_-{1_{\cal K}} &            &  {\cal K} \ar[ur]_-{H} &     
}}$\hskip 0.5cm
и

$\vcenter
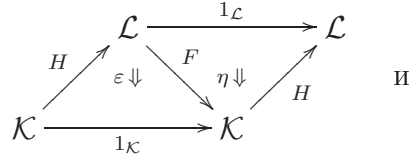
{\xymatrix{
          &  {\cal K} \ar[dr]^(0.44){H} \ar[rr]^-{1_{\cal K}} \ar@{}[d]|{\eta \, \Uparrow } &   \ar@{}[d]|{\varepsilon \, \Uparrow }    &   {\cal K}  \\
{\cal L} \ar[ur]^-{F} \ar[rr]_-{1_{\cal L}} &            &  {\cal L} \ar[ur]_-{F} &     
}}$\hskip 0.5cm
Они выражают по существу абстрактную 

\addvspace{0.1cm}
\noindent ситуацию сопряженности в 2-категории. 

\addvspace{0.1cm}
Имееется ввиду, что диаграммы коммутативны для 1- и 2- стрелок. Для 1-стрелок это очевидно. Для 2-стрелок это означает $(1_H*\varepsilon )\circ (\eta *1_H)=1_H$ и $(\varepsilon *1_F)\circ (1_F*\eta )=1_F$. Эти выражения называются обычно {\bf треугольными тождествами} и обозначаются короче $H\varepsilon \circ \eta H=1_H$ и $\varepsilon F\circ F\eta =1_F$.

\begin{proposition}
Для функторов $F:{\cal L}\to {\cal K}$ и $H:{\cal K}\to {\cal L}$ следующие утверждения эквивалентны:
\begin{itemize}
\item $F\dashv H$, 
\item существуют естественные преобразования $\varepsilon :FH\Rightarrow 1_{\cal K}$ и $\eta :1_{\cal L}\Rightarrow HF$, удовлетворяющие треугольным тождествам $H\varepsilon \circ \eta H=1_H$ и $\varepsilon F\circ F\eta =1_F$.
\end{itemize}
\end{proposition}
Доказательство.\\
(Из пункта 1 следует пункт 2): Рассмотрим диаграммы, выражающие универсальность стрелок $\varepsilon _Y$ и $\eta _X$\hskip 1cm
$\vcenter
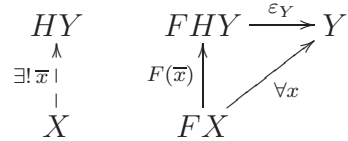
{\xymatrix{
HY  & FHY  \ar[r]^-{\varepsilon _{Y}} & Y \\
X \ar@{-->}[u]^-{\exists !\, \overline x}  &  FX \ar[ur]_-{\forall x} \ar[u]^-{F(\overline {x})} & 
}}$\\
$\vcenter
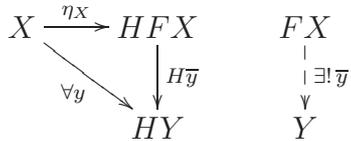
{\xymatrix{
X \ar[r]^-{\eta _X} \ar[dr]_-{\forall y} & HFX \ar[d]^-{H\overline y} & FX \ar@{-->}[d]^-{\exists !\, \overline y} \\
  & HY & Y
}}$\hskip 1cm Подставим в первую диаграмму $Y=FX$, $x=1_{FX}$, а во вторую $X=HY$, $y=1_{HY}$, получим требуемые равенства
$\varepsilon _{FX}\circ F\eta _X=1_{FX}$ и $H\varepsilon _Y\circ \eta _{HY}=1_{HY}$.

\addvspace{0.1cm}
\noindent (Из пункта 2 следует пункт 1): Определим отображения между $Hom$-множествами
$\xymatrix{{\cal K}(FX,Y)\ar@/^0.5pc/[r]^-{\theta _{X,Y}} & {\cal L}(X,HY) \ar@/^0.5pc/[l]^-{\theta ^*_{X,Y}} }$
следующим образом \\
$\left\{
\begin{array}{lr}
\theta _{X,Y}(f):=Hf\circ \eta _X, & f\in {\cal K}(FX,Y) \\
\theta ^*_{X,Y}(g):=\varepsilon _Y\circ Fg,  & g\in {\cal L}(X,HY)
\end{array}
\right. $ 

\addvspace{0.2cm}
\noindent Вычисляем $\theta ^*_{X,Y}(\theta _{X,Y}(f))=\varepsilon _Y\circ F(Hf\circ \eta _X)=\varepsilon _Y\circ FHf\circ F\eta _X=f\circ \varepsilon _{FX}\circ F\eta _X=f\circ 1_{FX}=f$ (мы использовали факт, что $\varepsilon $ естественное преобразование, и второе треугольное тождество).

\addvspace{0.1cm}
\noindent Аналогично, $\theta _{X,Y}(\theta ^*_{X,Y}(g))=H(\varepsilon _Y\circ Fg)\circ \eta _X=H\varepsilon _Y\circ HFg\circ \eta _X=H\varepsilon _{Y}\circ \eta _{HY}\circ g=1_{HY}\circ g=g$ (использовались естественность $\eta $ и первое треугольное тождество).

\addvspace{0.1cm}
\noindent Таким образом, $\theta _{X,Y}$ и $\theta ^*_{X,Y}$ взаимно обратные преобразования.

\addvspace{0.1cm}
\noindent Диаграмма $\vcenter{\xymatrix{
X \ar[d]_-{x} & Y  & {\cal K}(FX,Y) \ar[r]^-{\theta _{X,Y}} & {\cal L}(X,HY) \\
X' & Y' \ar[u]^-{y} & {\cal K}(FX',Y') \ar[u]^-{{\cal K}(Fx,y)} \ar[r]_-{\theta _{X',Y'}} & {\cal L}(X',HY') \ar[u]_-{{\cal L}(x,Hy)}
}}$
коммутативна, так как для всякой стрелки $h:FX'\to Y'$ выполняется $\theta _{X,Y}({\cal K}(Fx,y)(h))=\theta _{X,Y}(y\circ h\circ Fx)=H(y\circ h\circ Fx)\circ \eta _X=Hy\circ Hh\circ HFx\circ \eta _X=Hy\circ Hh\circ \eta _{X'}\circ x=Hy\circ \theta _{X',Y'}(h)\circ x={\cal L}(x,Hy)(\theta _{X',Y'}(h))$.

\addvspace{0.1cm}
\noindent Следовательно, $\theta _{X,Y}$ естественная биекция.      \hfill $\square $

\begin{definition}
Два объекта $A$ и $B$ в 2-категории ${\cal K}$ называются {\bf сопряженными}, если
\begin{itemize}
\item существуют 1-стрелки $\xymatrix{A \ar@/^0.5pc/[r]^-{f} & B \ar@/^0.5pc/[l]^-{g} }$,
\item существуют 2-стрелки $\varepsilon :gf\Rightarrow 1_A$, $\eta :1_B\Rightarrow fg$,
\item выполняются треуголные тождества $f\varepsilon \circ \eta f=1_f$, $\varepsilon g\circ g\eta =1_g$.
\end{itemize}

\noindent Стрелка $g$ называется {\bf левой сопряженной} к $f$, а $f$ {\bf правой сопряженной} к $g$, что обозначается $g\dashv f$. $\varepsilon $ называется {\bf коединицей}, а $\eta $ -- {\bf единицей} сопряжения.
\end{definition}

\begin{proposition}
2-функтор сохраняет ситуацию сопряжения.
\end{proposition}
Доказательство оставляется читателю.

\addvspace{0.2cm}
{\bf Примеры.}
\begin{enumerate}
\item В $2\text{-}\mathbf {Cat}$ ситуация сопряжения согласно определению 14 это обычная ситуация сопряжения категорий. 
\item В $2\text{-}\mathbf {Top}$ ситуация сопряжения это специальный случай слабой эквивалентности, так как в $2\text{-}\mathbf {Top}$ все 2-стрелки обратимы (то есть $\varepsilon $ и $\eta $ -- 2-изоморфизмы). Дополнительным соотношением на гомотопии будут треугольные тождества. Это находится в контрасте с $2\text{-}\mathbf {Cat}$, где сопряженность является более широким свойством, чем эквивалентность категорий. 
\end{enumerate}

\begin{assignment}
\begin{enumerate}
\item Разобраться в типах стрелок и корректности всех произведений, участвующих в треугольных тождествах. 
\item Рассмотреть, как 2-функтор действует на ситуацию сопряжения, во что отображает сопряженные стрелки, единицу и коединицу. Доказать утверждение 14.
\end{enumerate}
\end{assignment}

\chapter{Тензорные категории}

В третьей главе рассматриваются категории с дополнительной структурой (слабого) моноида на категории. Они имеют широкое применение для решения разных вопросов от гомологической алгебры и топологии до линейной логики и квантовой теории поля.

\section{Определение тензорной категории}

Напомним, что для категории ${\cal K}$ категория ${\cal K}\times {\cal K}$ состоит из пар объектов $(A,B)\in ({\cal K}\times {\cal K})_0$ и пар стрелок $(f,g)\in ({\cal K}\times {\cal K})_1$, где $A,B\in {\cal K}_0$, $f,g\in {\cal K}_1$. Единицами являются пары единиц $(1_A,1_B)$, композиция определена покомпонентно $(f,g)\circ (f',g'):=(f\circ f',g\circ g')$. 

\begin{definition}
Слабым {\bf моноидом на категории} ${\cal K}\in 2\text{-}\mathbf {Cat}_0$, называется сама категория ${\cal K}$, снабженная следующими стрелками:
\begin{itemize}
\item функтором $\otimes :{\cal K}\times {\cal K}\to {\cal K}$, называемым {\bf тензорным произведением}, 
\item функтором $I:\mathbf 1\to {\cal K}$, где $\mathbf 1$ -- {\bf единичная} категория, состоящая из одного объекта и одной стрелки (функтор $I$ выделяет один специальный объект $I$ в категории ${\cal K}$, который называется {\bf единичным}),
\item естественным изоморфизмом $\alpha _{A,B,C}:(A\otimes B)\otimes C\stackrel {\simeq }{\to }A\otimes (B\otimes C)$, называемым {\bf изоморфизмом ассоциативности},
\item естественным изоморфизмом $\lambda _A:I\otimes A\stackrel {\simeq }{\to }A$, называемым {\bf левым единичным изоморфизмом}, 
\item естественным изоморфизмом $\rho _A:A\otimes I\stackrel {\simeq }{\to }A$, называемым {\bf правым единичным изоморфизмом},
\end{itemize}
которые удовлетворяют некоторым условиям согласованности, называемым {\bf условиями когерентности}, и выражаемым двумя коммутативными диаграммами

\addvspace{0.25cm} 
$\vcenter{\xymatrix{
 & (A\otimes (B\otimes C))\otimes D \ar[dr]^-{\alpha _{A,B\otimes C,D}}   &   \\
((A\otimes B)\otimes C)\otimes D \ar[ur]^-{\alpha _{A,B,C}\otimes 1_D} \ar[d]_-{\alpha _{A\otimes B,C,D}} &    & A\otimes ((B\otimes C)\otimes D) \ar[d]^-{1_A\otimes \alpha _{B,C,D}}  \\
 (A\otimes B)\otimes (C\otimes D) \ar[rr]_-{\alpha _{A,B,C\otimes D}} &    &  A\otimes (B\otimes (C\otimes D))
}}$

\addvspace{0.25cm}
и \hskip 2cm $\vcenter{\xymatrix{
(A\otimes I)\otimes B \ar[rr]^-{\alpha _{A,I,B}} \ar[dr]_-{\rho _A\otimes 1_B} & & A\otimes (I\otimes B) \ar[dl]^-{1_A\otimes \lambda _B} \\
  & A\otimes B & 
}}$
\end{definition}

\addvspace{0.1cm}
\noindent Категория, на которой определена структура (слабого) моноида называется {\bf тензорной} или {\bf моноидальной}. 

\addvspace{0.1cm}
\noindent Несмотря на длину определения, в нем всего лишь вводятся правила оперирования объектами и стрелками в тензорной категории. Так, например, два объекта $A$ и $B$ можно перемножить $A\otimes B$, две стрелки $f:A\to A'$ и $g:B\to B'$ можно перемножить $f\otimes g:A\otimes B\to A'\otimes B'$, умножение удовлетворяет ассоциативному закону с точностью до эквивалентности, на единицу $I$ можно сокращать слева и справа с точностью до эквивалентности, сами эквивалентности, с точностью до которых выполняются операции, удовлетворяют некоторым естественным тождествам (условиям когерентности). На самом деле даже имеет место {\bf теорема когерентности} Маклейна-Сташефа, утверждающая, что в {\it тензорной категории все диаграммы, построенные при помощи стрелок $\lambda $, $\rho $, $\alpha $, единичных стрелок и операции тензорного умножения $\otimes $, коммутативны}.

Более непосредственным определением структуры (слабого) моноида на категории ${\cal K}$ было бы использование диаграмм (ассоциативности и единицы) 

\addvspace{0.2cm}
\noindent $\vcenter{\xymatrix{
{\cal K}\times ({\cal K}\times {\cal K}) \ar[r]^-{\simeq } \ar[d]_-{1_{{\cal K}}\times \otimes } & ({\cal K}\times {\cal K})\times {\cal K} \ar[r]^-{\otimes \times 1_{{\cal K}}} \ar@{}[d]|{\stackrel {\displaystyle \alpha }{\displaystyle \Longleftarrow }} & {\cal K}\times {\cal K} \ar[d]^-{\otimes } \\
{\cal K}\times {\cal K} \ar[rr]_-{\otimes } && {\cal K}
}}$
$\vcenter{\xymatrix{
{\cal K}\times {\cal K} \ar[r]^-{\otimes } & {\cal K} & {\cal K}\otimes {\cal K} \ar[l]_-{\otimes } \\
\mathbf 1\times {\cal K} \ar[u]^-{I\times 1_{\cal K}}_(0.65){\hskip 0.2cm\stackrel {\displaystyle \lambda }{\Rightarrow }} \ar[ur]_-{pr_2} & & {\cal K}\times \mathbf 1 \ar[u]_-{1_{\cal K}\times I}^(0.65){\stackrel {\displaystyle \rho }{\displaystyle \Leftarrow }\hskip 0.2cm} \ar[ul]^-{pr_1}
}}$

\noindent выполняющихся с точностью до естественных изоморфизмов $\alpha $, $\lambda $ и $\rho $ (которые в свою очередь удовлетворяют условиям когерентности). Читателю рекомендуется сравнить эти два определения и убедиться, что они полностью одинаковы.

\addvspace{0.2cm}
{\bf Примеры.}
\begin{enumerate}
\item Категория с бинарным произведением $\times $ и конечным объектом $\mathbf 1$ является тензорной с тензорным произведением $\otimes :=\times $ и единичным объектом $I:=\mathbf 1$.
\item Категория с бинарным копроизведением $\coprod $ и начальным объектом $\mathbf 0$ является тензорной с тензорным произведением $\otimes :=\coprod $ и единичным объектом $I:=\mathbf 0$.
\item Категория {\bf эндофункторов} $\mathbf {End}({\cal K}):=2\text{-}\mathbf {Cat}({\cal K},{\cal K})$ (то есть функторов, начинающихся и кончающихся в одной категории) является тензорной с тензорным произведением $\otimes :=\circ $ (композицией функторов) и $I:=1_{{\cal K}}$ (тождественным функтором). Если в двух первых примерах категории были {\it слабые}, то есть $\alpha $, $\lambda $ и $\rho $ были не тождественными естественными преобразованиями, то категория эндофункторов $\mathbf {End}({\cal K})$ является {\it строгой} (все преобразования $\alpha $, $\lambda $ и $\rho $ тождественные). На основе этой категории вводятся понятия монады и комонады.
\item Рассмотрим категорию векторных пространств $\mathbf {Vect}$ (над полем $k$). Построим категорию, объектами которой являются билинейные отображения $f:X_1\times X_2\to Y$, где $X_1,X_2, Y\in \mathbf {Vect}_0$ (векторные пространства), а стрелки $h:f\to f'$ -- линейные отображения $h:Y\to Y'$ такие, что диаграмма $\vcenter{\xymatrix{
X_1\times X_2 \ar[r]^-{f} \ar[dr]_-{f'} & Y \ar[d]^-{h} \\
  &  Y' 
}}$ 
коммутативна. (Билинейное отображение это отображение, линейное по каждому аргументу, например, $f(x_1,ax_2+bx'_2)=af(x_1,x_2)+bf(x_1,x'_2)$). В этой категории существует начальный объект, обозначаемый $\otimes :X_1\times X_2\to X_1\otimes X_2$ и называемый {\bf тензорным произведением} векторных пространств.

\addvspace{0.1cm}
\noindent $X_1\otimes X_2$ это векторное пространство, элементами которого являются {\it конечные} линейные комбинации вида $a_1x_{1_1}\otimes x_{2_1}+\cdots +a_nx_{1_n}\otimes x_{2_n}$, где $a_1,\dots a_n\in k$, $x_{1_1}, \dots , x_{1_n}\in X_1$, $x_{2_1},\dots ,x_{2_n}\in X_2$. Операция $\otimes $ рассматривается как билинейная операция на элементах пространств $X_1$, $X_2$ со значениями в третьем пространстве $X_1\otimes X_2$.

\addvspace{0.1cm}
\noindent Операция $\otimes :\mathbf {Vect}\times \mathbf {Vect}\to \mathbf {Vect}$ продолжается на стрелки так, что если $f:X_1\to X'_1$ и $g:X_2\to X'_2$ линейные отображения, то $f\otimes g:X_1\otimes X_2\to X'_1\otimes X'_2:a_1x_{1_1}\otimes x_{2_1}+\cdots +a_nx_{1_n}\otimes x_{2_n}\mapsto a_1f(x_{1_1})\otimes g(x_{2_1})+\cdots +a_nf(x_{1_n})\otimes g(x_{2_n})$ -- линейное отображение, называемое {\bf тензорным произведением} линейных отображений $f$ и $g$. Таким образом, тензорное произведение $\otimes $ является функтором двух аргументов на категории векторных пространств. Имеют место естественные изоморфизмы $\alpha _{X,Y,Z}:(X\otimes Y)\otimes Z\stackrel {\simeq }{\to }X\otimes (Y\otimes Z):(x\otimes y)\otimes z\mapsto x\otimes (y\otimes z)$, $\lambda _X:k\otimes X\stackrel {\simeq }{\to }X:1\otimes x\mapsto x$ и $\rho _X:X\otimes k\to X:x\otimes 1\mapsto x$. Следовательно, категория векторных пространств является тензорной категорией с тензорным произведением, являющимся (классическим) тензорным произведением векторных пространств, и единицей $I:=k$. Функтор тензорного умножения $X\otimes (-):\mathbf {Vect}\to \mathbf {Vect}$ имеет правый сопряженный $\mathbf {Vect}(X,-):\mathbf {Vect}\to \mathbf {Vect}$. Такие категории называются {\bf замкнутыми} по аналогии с декартово-замкнутыми категориями (в которых $\otimes :=\times $). Теория замкнутых тензорных категорий имеет развитое формальное исчисление и важна в частности в гомологической алгебре и линейной логике.  
\end{enumerate}

В некоторых теориях требуется, чтобы была некоторая связь между объектами $X\otimes Y$ и $Y\otimes X$. В этих случаях вводят естественный изоморфизм $\gamma :\otimes \stackrel {\simeq }{\to }\otimes \circ \tau :{\cal K}\times {\cal K}\to {\cal K}$, где $\tau :{\cal K}\times {\cal K}\to {\cal K}\times {\cal K}$ -- функтор, меняющий местами первую и вторую компоненты пары объектов и стрелок. Другими словами, вводятся стрелки $\gamma _{X,Y}:X\otimes Y\stackrel {\simeq }{\to }Y\otimes X$, естественно зависящие от аргументов $X, Y$. Если произведение $X\otimes Y\stackrel {\gamma _{X,Y}}{\to }Y\otimes X\stackrel {\gamma _{Y,X}}{\to }X\otimes Y$ есть тождественное естественное преобразование $1_{X\otimes Y}$ для всех $X,Y\in {\cal K}_0$, то категория ${\cal K}$ называется {\bf симметричной}, в противном случае {\bf несимметричной}. Если преобразование $\gamma $ отсутсвует, то симметричность или несимметричность тензорной категории неопределена. Требуется, чтобы $\gamma $ само удовлетворяло условиям согласованности (когерентности) с $\alpha $, $\lambda $, $\rho $. Мы не будем их выписывать.   

\addvspace{0.1cm}
Тензорные категории формируют 2-категорию, объектами которой они являются. Естественно рассматривать не все функторы между ними и не все естественные преобразования, а только согласованные с тензорной структурой. Мы также не будем на этом останавливаться, однако читателю рекомендуется подумать, как бы могли выглядеть эти условия согласованности.

\begin{assignment}
\begin{enumerate}
\item Убедиться, что определение 15, структуры слабого моноида на категории, выражает в точности то же самое, что и определение через диаграммы ассоциативности и единицы.
\item Написать уравнения, выражающие свойства тензорного произведения $\otimes :{\cal K}\times {\cal K}\to {\cal K}$ быть функтором.
\item Проверить корректность определения тензорных категорий в примерах 1 и 2.
\item Уточнить вид линейных изоморфизмов $\alpha _{X,Y,Z}$, $\lambda _X$, $\rho _X$ в примере 4 на произвольных элементах пространств.
\item Доказать, что отображение $\tau :{\cal K}\times {\cal K}\to {\cal K}\times {\cal K}$, меняющее местами первую и вторую компоненты объектов и стрелок, является функторным изоморфизмом.   
\end{enumerate}
\end{assignment}

\section{(Ко)моноиды}

Обычно в математике {\bf моноидом} называется множество $M$ с одной бинарной операцией $\cdot :M\times M\to M$, которая предполагается ассоциативной, то есть для всех $a,b,c\in M$ $(a\cdot b)\cdot c=a\cdot (b\cdot c)$, и имеющей единицу $e$, то есть для любого $a\in M$ $e\cdot a=a=a\cdot e$. Если записать это определение через диаграммы, то получим две диаграммы ({\it ассоциативности} и {\it единицы})

\noindent \hskip -0.35cm$\vcenter{\xymatrix{
M\times (M\times M) \ar[r]^-{\simeq } \ar[d]_-{1_M\times \cdot } & (M\times M)\times M \ar[r]^-{\cdot \times 1_M} & M\times M \ar[d]^-{\cdot } \\
M\times M  \ar[rr]_-{\cdot } &    &  M
}}$ \hskip -0.35cm
$\vcenter{\xymatrix{
M\times M \ar[r]^-{\cdot } & M & M\times M \ar[l]_-{\cdot } \\
\mathbf 1\times M \ar[u]^-{e\times 1_M} \ar[ur]_-{pr_2} & & M\times \mathbf 1 \ar[u]_-{1_M\times e} \ar[ul]^-{pr_1}
}}$

\noindent где $\mathbf 1$ есть одноэлементное множество (конечный объект в категории множеств $\mathbf {Set}$).

По аналогии определяются моноиды в любой тензорной категории. Комоноиды получаются 'обращением стрелок'.
\begin{definition}
Пусть $({\cal K},\otimes ,I,\alpha ,\lambda ,\rho )$ -- тензорная категория.
\begin{itemize}
\item {\bf Моноидом} в ${\cal K}$ называется объект $M\in {\cal K}_0$ вместе с двумя стрелками: {\bf умножением}  $\mu :M\otimes M\to M$ и {\bf единицей} $\eta :I\to M$ такими, что выполняются уравнения, задаваемые диаграммами ассоциативности и единицы 

\noindent \hskip -1.5cm$\vcenter{\xymatrix{
M\otimes (M\otimes M) \ar[r]^-{\alpha ^{-1}_{M,M,M}} \ar[d]_-{1_M\otimes \mu } & (M\otimes M)\otimes M \ar[r]^-{\mu \otimes 1_M} & M\otimes M \ar[d]^-{\mu } \\
M\otimes M  \ar[rr]_-{\mu } &    &  M
}}$ \hskip -0.35cm
$\vcenter{\xymatrix{
M\otimes M \ar[r]^-{\mu } & M & M\otimes M \ar[l]_-{\mu } \\
I\otimes M \ar[u]^-{\eta \otimes 1_M} \ar[ur]_-{\lambda _M} & & M\otimes I\ar[u]_-{1_M\otimes \eta } \ar[ul]^-{\rho _M}
}}$
\item {\bf Комоноидом} в ${\cal K}$ называется объект $C\in {\cal K}_0$ вместе с двумя стрелками: {\bf коумножением}  $\delta :C\to C\otimes C$ и {\bf коединицей} $\varepsilon :C\to I$ такими, что выполняются уравнения, задаваемые диаграммами коассоциативности и коединицы 

\noindent \hskip -0.3cm$\vcenter{\xymatrix{
C\otimes (C\otimes C)  & (C\otimes C)\otimes C \ar[l]_-{\alpha _{C,C,C}} & C\otimes C \ar[l]_-{\delta \otimes 1_C}  \\
C\otimes C \ar[u]^-{1_C\otimes \delta }  &    &  C \ar[u]_-{\delta  } \ar[ll]^-{\delta  }
}}$ \hskip -0.35cm
$\vcenter{\xymatrix{
C\otimes C \ar[d]_-{\varepsilon \otimes 1_C} & C \ar[dl]^-{\lambda ^{-1}_C} \ar[dr]_-{\rho ^{-1}_C} \ar[l]_-{\delta } \ar[r]^-{\delta } & C\otimes C \ar[d]^-{1_C\otimes \varepsilon } \\
I\otimes C  & & C\otimes I 
}}$
\end{itemize}
\end{definition} 

\addvspace{0.2cm}
{\bf Примеры.}
\begin{enumerate}
\item Обычными (классическими) моноидами являются, например, множество натуральных чисел с операцией сложения и нулем в качестве единицы $(\Bbb N,+,0)$, множество натуральных чисел с операцией умножения и обычной единицей $(\Bbb N,\cdot ,1)$, множество квадратных матриц (над кольцом с единицей), множество преобразований множества $A$ в себя $(\mathbf {Set}(A,A),\circ ,1_A)$. Все это моноиды в тензорной категории $\mathbf {Set}$ с декартовым произведением $\times $ в качестве тензорного умножения $\otimes $ и единицей $I:=\mathbf 1$.  
\item Моноиды и комоноиды в категории эндофункторов \, $\mathbf {End}({\cal K}):=$ $2\text{-}\mathbf {Cat}({\cal K},{\cal K})$ \, являются соответственно монадами и комонадами (смотрите следующий параграф). Это пример моноидов и комоноидов, которые не являются множествами.
\item Единственно возможной структурой комоноида на множестве $A\in \mathbf {Set}_0$ ($\mathbf {Set}$ рассматривается как тензорная категория по отношению к декартову произведению $\times $ и единице $\mathbf 1$) является {\bf диагональное отображение} $\delta :A\to A\times A:a\mapsto (a,a)$ (коумножение) вместе с отображением в одноэлементное множество $!:A\to \mathbf 1$ (коединица). 
\item Рассмотрим векторное пространство $k[X]\in \mathbf {Vect}_0$, порожденное множеством $X$ (то есть $X$ является его базисом), в тензорной категории $(\mathbf {Vect},\otimes ,k)$. Оно допускает структуру комоноида с коумножением $\delta :k[X]\to k[X]\otimes k[X]:x\mapsto x\otimes x$ и коединицей $\varepsilon :X\to k:x\mapsto 1$ (отображения заданы на базисных элементах и распространяются на другие по линейности).
\end{enumerate}

\begin{assignment}
\begin{enumerate}
\item Доказать, что моноид в категории моноидов $\mathbf {Mon}$, рассматриваемой с декартовым произведением в качестве тензорного умножения и единицей -- единичным моноидом (состоящим из одного элемента), является коммутативным моноидом в $\mathbf {Set}$.
\item Проверить диаграммы коассоциативности и коединицы для структуры комоноида на множестве (пример 3). Доказать единственность структуры комоноида на множестве в тензорной категории $(\mathbf {Set},\times ,\mathbf 1)$ (используя диаграмму коединицы).
\item Проверить корректность задания структуры комоноида на векторном пространстве в примере 4.
\end{enumerate}
\end{assignment}

\section{(Ко)монады}

Фиксируем объект $A$ в 2-категории ${\cal K}$. По определению 2-категории $Hom$-множество ${\cal K}(A,A)$ представляет собой 1-категорию с объектами -- 1-стрелками в ${\cal K}$ типа $f:A\to A$ и стрелками -- 2-стрелками в ${\cal K}$ типа $\xymatrix{A \ar@<1.25ex>[r]^-{f}_-{\alpha \Downarrow } \ar@<-1.25ex>[r]_-{g}  &    A}$. Эта категория является строгой тензорной с тензорным произведением 
$\left\{
\begin{array}{lr}
f\otimes g:=f\circ g & f,g\in {\cal K}(A,A)_0 \\
\alpha \otimes \beta:=\alpha *\beta & \alpha ,\beta \in {\cal K}(A,A)_1 
\end{array}
\right. $ 
и единицей $I:=1_A$. Монады это моноиды, а комонады -- комоноиды в категории ${\cal K}(A,A)$. Замечательно, что в {\it ситуации сопряжения} $\xymatrix{A \ar@/^0.5pc/[r]^-{f}_-{\top }  &   B \ar@/^0.5pc/[l]^-{g} }$ стрелка $f\circ g:B\to B$ всегда является монадой в ${\cal K}(B,B)$, а стрелка $g\circ f:A\to A$ -- комонадой в ${\cal K}(A,A)$. Дадим точные определения.

\begin{definition}
Пусть ${\cal K}$ -- 2-категория. 1-стрелка $f:A\to A$ называется
\begin{itemize}
\item {\bf монадой}, если заданы две 2-стрелки $\mu :f\circ f\Rightarrow f$ (умножение) и $\eta :1_A\Rightarrow f$ (единица) такие, что диаграммы ассоциативности и единицы 
\hskip 0.35cm$\vcenter{\xymatrix{
f\circ f\circ f \ar[r]^-{1_f*\mu } \ar[d]_-{\mu *1_f} & f\circ f \ar[d]^-{\mu } \\
f\circ f \ar[r]_-{\mu } &  f
}}$
$\vcenter{\xymatrix{
1_A\circ f \ar[r]^-{\eta *1_f} \ar[dr]_-{=} & f\circ f \ar[d]^-{\mu } & f\circ 1_A \ar[l]_-{1_f*\eta } \ar[dl]^-{=} \\
   & f & 
}}$\hskip 0.35cm
коммутативны.
\item {\bf комонадой}, если заданы две 2-стрелки $\delta :f\Rightarrow f\circ f$ (коумножение) и $\varepsilon :f\Rightarrow 1_A$ (коединица) такие, что диаграммы коассоциативности и коединицы 

\noindent $\vcenter{\xymatrix{
f\circ f\circ f  & f\circ f \ar[l]_-{1_f*\delta } \\
f\circ f \ar[u]^-{\delta *1_f} &  f \ar[l]^-{\delta } \ar[u]_-{\delta }
}}$
$\vcenter{\xymatrix{
1_A\circ f   & f\circ f \ar[l]_-{\varepsilon *1_f} \ar[r]^-{1_f*\varepsilon } & f\circ 1_A   \\
   & f \ar[ul]^-{=} \ar[u]_-{\delta } \ar[ur]_-{=} & 
}}$\hskip 0.35cm
коммутативны.
\end{itemize} 
\end{definition}

\begin{proposition}
Пусть $\xymatrix{A \ar@/^0.5pc/[r]^-{f}_-{\top }  &   B \ar@/^0.5pc/[l]^-{g} }$ -- ситуация сопряжения (с единицей $\eta :1_B\Rightarrow f\circ g$ и коединицей $\varepsilon :g\circ f\Rightarrow 1_A$) в 2-категории ${\cal K}$, тогда
\begin{itemize}
\item $(f\circ g,1_f*\varepsilon *1_g,\eta )$ является монадой на стрелке $f\circ g:B\to B$ с умножением $\mu :=1_f*\varepsilon *1_g:f\circ g\circ f\circ g\Rightarrow f\circ g$ и единицей $\eta :1_B\Rightarrow f\circ g$,
\item $(g\circ f,1_g*\eta *1_f,\varepsilon )$ является комонадой на стрелке $g\circ f:A\to A$ с коумножением $\delta :=1_g*\eta *1_f:g\circ f\Rightarrow g\circ f\circ g\circ f$ и коединицей $\varepsilon :g\circ f\Rightarrow 1_A$.
\end{itemize}
\end{proposition}
Доказательство.
Докажем первую часть утверждения.

\addvspace{0.15cm}
\noindent (диаграмма ассоциативности): Нужно доказать $\mu \circ (1_{f\circ g}*\mu )=\mu \circ (\mu *1_{f\circ g})$.
Подставим в левую и правую части $\mu $ и вычислим выражения, используя правила алгебры в 2-категории. \\
Левая часть: $(1_f*\varepsilon *1_g)\circ (1_{f\circ g}*1_f*\varepsilon *1_g)=1_f*(\varepsilon \circ (1_{g\circ f}*\varepsilon))*1_g=1_f*\varepsilon *\varepsilon *1_g$.\\
Правая часть: $(1_f*\varepsilon *1_g)\circ (1_f*\varepsilon *1_g*1_{f\circ g})=1_f*(\varepsilon \circ (\varepsilon *1_{g\circ f}))*1_g=1_f*\varepsilon *\varepsilon *1_g$. [В вычислениях использовалось в частности, что $1_{f\circ g}=1_f*1_g$ и $\varepsilon \circ (\varepsilon *1_{g\circ f})=(1_{1_A}\circ \varepsilon )*(\varepsilon \circ 1_{g\circ f})=\varepsilon *\varepsilon =(\varepsilon \circ 1_{g\circ f})*(1_{1_A}\circ \varepsilon )=\varepsilon \circ (1_{g\circ f}*\varepsilon )$]

\addvspace{0.2cm}
\noindent (диаграмма единицы): Нужно доказать $\mu \circ (\eta *1_{f\circ g})=1_{f\circ g}$ и $\mu \circ (1_{f\circ g}*\eta )=1_{f\circ g}$.\\
Подставляем $\mu $ и вычисляем: $(1_f*\varepsilon *1_g)\circ (\eta *1_{f\circ g})=((1_f*\varepsilon )\circ (\eta *1_f))*1_g=1_f*1_g=1_{f\circ g}$ и $(1_f*\varepsilon *1_g)\circ (1_{f\circ g}*\eta )=1_f*((\varepsilon *1_g)\circ (1_g*\eta ))=1_f*1_g=1_{f\circ g}$. [Использовались треугольные тождества]        \hfill $\square $  

\addvspace{0.2cm}
В категории $2\text{-}\mathbf {Cat}$ верно также и обратное утверждение, что всякая монада или комонада получаются из ситуации сопряжения, причем существует много ситуаций сопряжения, дающих данную монаду или комонаду. Пусть, например, $(T,\mu ,\eta )$ монада на категории ${\cal K}$. Всегда существует непустая категория сопряжений, дающих данную монаду в смысле утверждения 15. Объектами такой категории служат сопряжения $\vcenter{\xymatrix{
{\cal L} \ar@/^0.5pc/[d]^-{U} \\
{\cal K} \ar@/^0.5pc/[u]^-{F}_-{\dashv }
}}$
а стрелками функторы $\Phi :{\cal L}\to {L'}$, для которых треугольники в диаграмме $\vcenter{\xymatrix{
{\cal L} \ar[rr]^-{\Phi } \ar@/^0.5pc/[dr]^(0.75){U} & & {\cal L'} \ar@/^0.5pc/[dl]^(0.2){U'} \\
    & {\cal K} \ar@/^0.5pc/[ul]^(0.65){F}_(0.38){\hskip -0.08cm\dashv } \ar@/^0.5pc/[ur]^(0.4){F'\hskip -0.1cm}_(0.7){\hskip -0.1cm\dashv } &  
}}$
коммутативны. Имеются {\it начальный} и {\it конечный объекты} в этой категории, так называемые, {\bf категория Клейсли} и {\bf категория алгебр Эйленберга-Мура}. Мы не будем на этом останавливаться.  

\addvspace{0.2cm}
{\bf Примеры.}
\begin{enumerate}
\item Пусть $\xymatrix{\mathbf {Set} \ar@/^0.5pc/[r]^-{F}_-{\bot } & \mathbf {Mon} \ar@/^0.5pc/[l]^-{U} }$ ситуация сопряжения между категориями множеств и моноидов. Тогда $(UF,U\varepsilon F,\eta )$ -- монада на категории $\mathbf {Set}$. Умножение $U\varepsilon F(A):UFUF(A)\to UF(A)$ сопоставляет слову $(a^1_1\cdots a^1_{n_1})\cdots (a^m_1\cdots a^m_{n_m})$ в $UFUF(A)$ слово $a^1_1\cdots a^1_{n_1}\cdots a^m_1\cdots a^m_{n_m}$ в $UF(A)$, где все $a^i_j$ принадлежат множеству $A$. Единица $\eta _A:A\to UF(A)$ -- единица сопряжения, которая сопоставляет элементу $a\in A$ слово $a\in UF(A)$.\\
$(FU,F\eta U,\varepsilon )$ -- комонада на категории $\mathbf {Mon}$. Коумножение $F\eta U(M):FU(M)\to FUFU(M)$ сопоставляет слову $m_1\cdots m_n$ в $FU(M)$ слово $(m_1)\cdots (m_n)$ в $FUFU(M)$, где $m_1,\dots ,m_n\in M$. Коединица $\varepsilon _M:FU(M)\to M$ -- коединица сопряжения, которая сопоставляет слову $m_1\cdots m_n$ в $FU(M)$ элемент моноида $M$, являющийся произведением элементов $m_1,\dots ,m_n$ в $M$.    
\item Пусть ${\cal P}:\mathbf {Set}\to \mathbf {Set}$ -- {\it ковариантный} функтор {\it множество-степень} (то есть если $f:A\to B$ -- функция, то ${\cal P}(f):{\cal P}(A)\to {\cal P}(B)$ -- функция, назначающая подмножеству $X\subset A$ его образ $f(X)\subset B$). Тогда $({\cal P},\bigcup ,\eta )$ -- монада на категории множеств $\mathbf {Set}$. Умножение $\bigcup :{\cal P}({\cal P}(A))\to {\cal P}(A)$ сопоставляет множеству подмножеств $A$ их объединение. Единица $\eta _A:A\to {\cal P}(A)$ сопоставляет элементу $a\in A$ одноэлементное множество $\{a\}\in {\cal P}(A)$.
\end{enumerate}

Монады на категории $\mathbf {Set}$ представляют собой обобщенные {\it алгебры} (множества с операциями), которые, как видно из примера 2, могут иметь бесконечноместные операции. На самом деле {\it любой} алгебраической системе (такой как $\mathbf {Mon}$, $\mathbf {Grp}$, $\mathbf {Rng}$ и так далее) соответствует определенная монада на категории множеств $\mathbf {Set}$. Произвольной монаде на категории множеств может не соответствовать никакая классическая алгебраическая система, но это будет алгебра в некотором смысле. Смысл комонад менее очевиден, однако известны их применения в топологии и компьютерных науках.

\begin{assignment}
\begin{enumerate}
\item Доказать вторую часть утверждения 15.
\item Проверить выполнение аксиом монады для ковариантного функтора множество-степень ${\cal P}:\mathbf {Set}\to \mathbf {Set}$ (пример 2).
\end{enumerate}
\end{assignment}

\section{(Ко)алгебры}

Алгебры являются моноидами, а коалгебры комоноидами в категории векторных пространств $\mathbf {Vect}$.

\begin{definition}
Пусть $(\mathbf {Vect},\otimes ,k)$ категория векторных пространств над полем $k$, рассматриваемая с тензорным произведением $\otimes $ и единицей $I:=k$. Тогда
\begin{itemize}
\item {\bf алгебра} $A$ это объект в $\mathbf {Vect}$, являющийся моноидом, то есть векторное пространство $A$, снабженное дополнительной структурой: 
\begin{itemize}
\item стрелкой $\mu :A\otimes A\to A$ (умножением)
\item и стрелкой $\eta :k\to A$ (единицей) такими, что следующие диаграммы коммутативны

\addvspace{0.2cm} 
\hskip -2cm$\vcenter{\xymatrix{
(A\otimes A)\otimes A \ar[r]^-{\alpha _{A,A,A}}_-{\simeq } \ar[d]_-{\mu \otimes 1_A} & A\otimes (A\otimes A) \ar[r]^-{1_A\otimes \mu } & A\otimes A \ar[d]^-{\mu } \\
A\otimes A \ar[rr]_-{\mu } & & A
}}$
\hskip 0.35cm$\vcenter{\xymatrix{
k\otimes A \ar[r]^-{\eta \otimes 1_A} \ar[dr]_-{\lambda _A}^-{\simeq } & A\otimes A \ar[d]^-{\mu } &  A\otimes k \ar[l]_-{1_A\otimes \eta } \ar[dl]^-{\rho _A}_-{\simeq } \\
  &  A  &  
}}$

\end{itemize}  
\item {\bf коалгебра} $C$ это объект в $\mathbf {Vect}$, являющийся комоноидом, то есть векторное пространство $C$, снабженное дополнительной структурой: 
\begin{itemize}
\item стрелкой $\delta :C\to C\otimes C$ (коумножением)
\item и стрелкой $\varepsilon :A\to k$ (коединицей) такими, что следующие диаграммы коммутативны

\addvspace{0.2cm} 
\hskip -2cm$\vcenter{\xymatrix{
(C\otimes C)\otimes C   & C\otimes (C\otimes C) \ar[l]_-{\alpha ^{-1}_{C,C,C}}^-{\simeq }  & C\otimes C  \ar[l]_-{1_C\otimes \delta } \\
C\otimes C  \ar[u]^-{\delta \otimes 1_C} & & C \ar[u]_-{\delta } \ar[ll]^-{\delta }
}}$
\hskip 0.35cm$\vcenter{\xymatrix{
k\otimes C   & C\otimes C \ar[l]_-{\varepsilon \otimes 1_C} \ar[r]^-{1_C\otimes \varepsilon } &  C\otimes k   \\
  &  C \ar[u]_-{\delta } \ar[ul]^-{\lambda ^{-1}_C}_-{\simeq } \ar[ur]_-{\rho ^{-1}_C}^-{\simeq } &  
}}$

\end{itemize}  
\end{itemize}
\end{definition}

\addvspace{0.2cm}
{\bf Примеры.}
\begin{enumerate}
\item Одномерное векторное пространство $k$ является алгеброй и коалгеброй с умножением $\mu :k\otimes k\to k:r\otimes s\mapsto rs$, единицей $\eta :k\to k:r\mapsto r$, коумножением $\delta :k\to k\otimes k:r\cdot 1\mapsto r(1\otimes 1)$ и коединицей $\varepsilon :k\to k: r\mapsto r$.
\item Пусть $(X):=\{1,X, X^2,\dots ,X^n,\dots\}$ {\it моноид с одной образующей} $X$ (и умножением $X^n\cdot X^m:=X^{n+m}$), пусть $k[(X)]$ векторное пространство с элементами моноида $(X)$  в качестве базиса (например, $r_0+r_1X+r_2X^2+r_5X^5$ является элементом $k[(X)]$). Тогда $k[(X)]$ является алгеброй и коалгеброй с умножением $\mu :k[(X)]\otimes k[(X)]\to k[(X)]:X^n\otimes X^m\mapsto X^{n+m}$, единицей $\eta :k\to k[(X)]:1\mapsto 1$, коумножением $\delta :k[(X)]\to k[(X)]\otimes k[(X)]:X^n\mapsto \sum \limits _{\stackrel {\scriptstyle i_1+i_2=n}{i_1,i_2\ge 0}} X^{i_1}\otimes X^{i_2}$ и коединицей $\varepsilon :k[(X)]\to k:\left\{
\begin{array}{lr}
X^n\mapsto 0 & \text{если } n>0\\
X^n\mapsto 1 & \text{если } n=0
\end{array}
\right. $ (все опреции заданы на базисных элементах и распространяются на другие по линейности).
\item Пусть $A\in \mathbf{Vect}_0$ векторное пространство, $A^{\otimes 2}:=A\otimes A$ тензорный квадрат $A$, то есть универсальное решение задачи для билинейных отображений $\vcenter{\xymatrix{
A\times A \ar[r]^-{\otimes } \ar[dr]_-{\forall f} & A\otimes A \ar@{-->}[d]^-{\exists !\, h} \\
   &  B
}}$,
$A^{\otimes 3}:=A\otimes A\otimes A$ тензорный куб $A$, то есть универсальное решение задачи для трилинейных отображений $\vcenter{\xymatrix{
A\times A\times A \ar[r]^-{\otimes } \ar[dr]_-{\forall f} & A\otimes A\otimes A \ar@{-->}[d]^-{\exists !\, h} \\
   &  B
}}$ (элементы $A^{\otimes 3}$ это конечные суммы вида $r_1a_{1_1}\otimes a_{1_2}\otimes a_{1_3}+\cdots +r_na_{n_1}\otimes a_{n_2}\otimes a_{n_3}$), \dots .
Рассмотрим копроизведение всех тензорных степеней пространства $A$: $T(A):=\coprod \limits _{0\le i< \infty }A^{\otimes i}:=\bigoplus \limits _{0\le i< \infty }A^{\otimes i}$. Элементы $T(A)$ -- конечные суммы {\it тензоров} разных степеней с коэффициентами из поля $k$ (например, $r_0+r_3a_1\otimes a_2\times a_3+r_4\bar a_1\otimes \bar a_2\otimes \bar a_3\otimes \bar a_4$ является элементом $T(A)$). Векторное пространство $T(A)$ имеет структуру алгебры и коалгебры с умножением $\mu :T(A)\otimes T(A)\to T(A):(a_1\otimes \cdots \otimes a_n)\otimes (\bar a_1\otimes \cdots \otimes \bar a_m)\mapsto a_1\otimes \cdots \otimes a_n\otimes \bar a_1\otimes \cdots \otimes \bar a_m$, единицей $\eta :k\to T(A):1\mapsto 1$, коумножением $\delta :T(A)\to T(A)\otimes T(A):a_1\otimes \cdots \otimes a_n\mapsto (1)\otimes (a_1\otimes \cdots \otimes a_n)+(a_1)\otimes (a_2\otimes \cdots \otimes a_n)+(a_1\otimes a_2)\otimes (a_3\otimes \cdots \otimes a_n)+\cdots +(a_1\otimes \cdots \otimes a_n)\otimes (1)$ и коединицей $\varepsilon :T(A)\to k:\left\{\begin{array}{lr}
1\mapsto 1 &  \\
a_1\otimes \cdots \otimes a_n\mapsto 0 & \text{для } n>0
\end{array}
\right.$ (операции определены на мономах или на базисных элементах и распространяются на другие элементы по линейности). Имеется другая (более сложная) структура коалгебры на $T(A)$, согласованная со структурой алгебры.
\end{enumerate}

\begin{definition}
{\bf Тензорным произведением}
\begin{itemize}
\item {\bf алгебр} $(A_1,\mu _1,\eta _1)$ и $(A_2,\mu _2,\eta _2)$ называется алгебра $(A_1\otimes A_2,\mu ,\eta )$, в которой 
\begin{itemize}
\item $A_1\otimes A_2$ -- тензорное произведение векторных пространств, 
\item $\mu $ -- умножение, определенное как композиция\\ 
$\vcenter{\xymatrix{A_1\otimes A_2\otimes A_1\otimes A_2 \ar[r]^-{\tau _{23}} & A_1\otimes A_1\otimes A_2\otimes A_2 \ar[r]^-{\mu _1\otimes \mu _2} & A_1\otimes A_2}}$, где\\ 
$\tau _{23}:A_1\otimes A_2\otimes A_1\otimes A_2\to A_1\otimes A_1\otimes A_2\otimes A_2:a_1\otimes a_2\otimes \bar a_1\otimes \bar a_2\mapsto a_1\otimes \bar a_1\otimes a_2\otimes \bar a_2$ -- канонический изоморфизм, меняющий местами второй и третий сомножители в тензорном произведении мономов,
\item $\eta $ -- единица, определенная как композиция\\
$\vcenter{\xymatrix{k \ar[r]^-{\delta _k} & k\otimes k \ar[r]^-{\eta _1\otimes \eta _2} & A_1\otimes A_2 }}$, где $\delta _k:k\to k\otimes k:1\mapsto 1\otimes 1$.
\end{itemize}
\item {\bf коалгебр} $(C_1,\delta _1,\varepsilon _1)$ и $(C_2,\delta  _2,\varepsilon _2)$ называется коалгебра $(C_1\otimes C_2,\delta  ,\varepsilon )$, в которой 
\begin{itemize}
\item $C_1\otimes C_2$ -- тензорное произведение векторных пространств, 
\item $\delta $ -- коумножение, определенное как композиция\\ 
$\vcenter{\xymatrix{
C_1\otimes C_2 \ar[r]^-{\delta _1\otimes \delta _2} & C_1\otimes C_1\otimes C_2\otimes C_2 \ar[r]^-{\tau _{23}} & C_1\otimes C_2\otimes C_1\otimes C_2
}}$, где\\ 
$\tau _{23}:C_1\otimes C_2\otimes C_1\otimes C_2\to C_1\otimes C_1\otimes C_2\otimes C_2:c_1\otimes c_2\otimes \bar c_1\otimes \bar c_2\mapsto c_1\otimes \bar c_1\otimes c_2\otimes \bar c_2$ -- канонический изоморфизм, меняющий местами второй и третий сомножители в тензорном произведении мономов,
\item $\varepsilon $ -- коединица, определенная как композиция\\
$\vcenter{\xymatrix{
C_1\otimes C_2 \ar[r]^-{\varepsilon _1\otimes \varepsilon _2} & k\otimes k \ar[r]^-{\mu _k} & k
}}$, где $\mu _k:k\otimes k\to k:1\otimes 1\mapsto 1$.
\end{itemize}
\end{itemize}
\end{definition}

Алгебры и коалгебры формируют тензорные подкатегории категории векторных пространств $\mathbf {Vect}$ с тензорным произведением, определенным выше, и единицей $k$. Стрелками в этих категориях выступают {\it гомоморфизмы} алгебр и коалгебр.

\begin{definition}
\begin{itemize}
\item {\bf Гомоморфизмом алгебр} $f:(A_1,\mu _1,\eta _1)\to (A_2,\mu _2,\eta _2)$ называется линейное отображение $f:A_1\to A_2$, коммутирующее с операциями умножения и единицы, то есть такое, для которого следующие диаграммы коммутативны \\
\vskip 0.0cm \hskip 1cm$\vcenter{\xymatrix{
A_1\otimes A_1 \ar[r]^-{f\otimes f} \ar[d]_-{\mu _1} & A_2\otimes A_2 \ar[d]^-{\mu _2} \\ 
A_1  \ar[r]_-{f}          &    A_2
}}$\hskip 1cm
$\vcenter{\xymatrix{
A_1 \ar[rr]^-{f} & & A_2 \\
  & k \ar[ul]^-{\eta _1} \ar[ur]_-{\eta _2} & 
}}$
\item {\bf Гомоморфизмом коалгебр} $f:(C_1,\delta _1,\varepsilon _1)\to (C_2,\delta _2,\varepsilon _2)$ называется линейное отображение $f:C_1\to C_2$, коммутирующее с операциями коумножения и коединицы, то есть такое, для которого следующие диаграммы коммутативны  \\
\vskip 0.0cm \hskip 1cm$\vcenter{\xymatrix{
C_1\otimes C_1 \ar[r]^-{f\otimes f}  & C_2\otimes C_2  \\ 
C_1  \ar[r]_-{f} \ar[u]^-{\delta _1}         &    A_2 \ar[u]_-{\delta _2}
}}$\hskip 1cm
$\vcenter{\xymatrix{
A_1 \ar[rr]^-{f} \ar[dr]_-{\varepsilon _1} & & A_2 \ar[dl]^-{\varepsilon _2} \\
  & k  & 
}}$
\end{itemize}
\end{definition}

\begin{assignment}
\begin{enumerate}
\item Убедиться, что определение алгебры, данное в тексте, совпадает с обычным определением алгебры над полем $k$ (смотрите учебник алгебры).
\item Проверить корректность определения структуры алгебры и коалгебры на пространстве $k$ (пример 1).
\item Проверить выполнение аксиом алгебры и коалгебры в примере 2.
\item Проверить корректность определения тензорного произведения алгебр и коалгебр (определение 19).
\end{enumerate}
\end{assignment}

\section{Алгебры Хопфа}

Алгебры Хопфа обычно определяются в тензорной категории векторных пространств $(\mathbf {Vect},\otimes ,k)$, но также имеют смысл в произвольной симметричной тензорной категории. Мы будем следовать второму подходу, чтобы была видна преемственность алгебр Хопфа от обычных групп. Алгебры Хопфа называются еще {\it квантовыми группами} и аналогично группам часто возникают как симметрии (некоммутативных) пространств.

Моноиды и комоноиды в симметричной тензорной категории $({\cal K},\otimes ,I)$ можно тензорно умножать точно так же как это было введено в предыдущем параграфе для алгебр и коалгебр ($\tau _{23}:=1_{A_1}\otimes \gamma _{A_2,A_1}\otimes 1_{A_2}$ для моноидов, и аналогично для комоноидов). Единичный объект $I\in {\cal K}_0$ имеет структуру моноида и комоноида с умножением $(\mu _I:I\otimes I\to I)=\lambda _I$, единицей $(\eta :I\to I)=1_I$, коумножением $(\delta :I\to I\otimes I)=\lambda ^{-1}_I$ и коединицей $(\varepsilon :I\to I)=1_I$. Доказательство этих фактов техническое, и мы его опустим (в основном два вспомогательных факта используются, что $\lambda _I=\rho _I$ и $1_I\otimes \lambda _I=(\lambda _I\otimes 1_I)\circ \alpha _{I,I,I}$).

\begin{definition}
Пусть $({\cal K},\otimes ,I)$ -- симметричная тензорная категория. 
{\bf Бимоноидом} в ${\cal K}$ называется объект $B\in {\cal K}_0$, являющийся одновременно моноидом и комоноидом, причем умножение $\mu :B\otimes B\to B$ и единица $\eta :I\to B$ являются гомоморфизмами комоноидов (и наоборот, что эквивалентно, коумножение $\delta :B\to B\otimes B$ и коединица $\varepsilon :B\to I$ являются гомоморфизмами моноидов), то есть следующие диаграммы коммутативны \hskip 0.5cm$\vcenter{\xymatrix{
B\otimes B \ar[r]^-{\mu } \ar[d]_-{\tau _{23}\circ (\delta \otimes \delta )} & B \ar[d]^-{\delta } \\
B\otimes B\otimes B\otimes B \ar[r]_-{\mu \otimes \mu } & B\otimes B
}}$\hskip 0.8cm
$\vcenter{\xymatrix{
B\otimes B \ar[rr]^-{\mu } \ar[dr]_-{\lambda _I\circ (\varepsilon \otimes \varepsilon )} && B \ar[dl]^-{\varepsilon } \\ 
 & I & 
}}$\\ 
\vskip 0.0cm \hskip4.5cm$\vcenter{\xymatrix{
I \ar[r]^-{\eta } \ar[d]_-{\lambda ^{-1}_I} & B \ar[d]^-{\delta } \\
I\otimes I \ar[r]_-{\eta \otimes \eta } & B\otimes B
}}$\hskip 1cm
$\vcenter{\xymatrix{
I \ar[rr]^-{\eta } \ar[dr]_-{1_I} && B \ar[dl]^-{\varepsilon } \\
 & I & 
}}$

\end{definition} 

Из определения следует, что бимоноиды являются моноидами в категории комоноидов (или, что то же самое, комоноидами в категории моноидов). Следовательно, категории моноидов на категории комоноидов и комоноидов на категории моноидов изоморфны.

\addvspace{0.1cm}
Пусть $C\in {\cal K}_0$ и $M\in {\cal K}_0$ соответственно комоноид и моноид в категории $({\cal K},\otimes ,I)$. На $Hom$-множестве ${\cal K}(C,M)$ имеется операция умножения $*$, называемая {\bf сверткой} или {\bf конволюцией}, определенная как $f*g:=\mu \circ (f\otimes g)\circ \delta $ или через диаграмму $C\stackrel {\delta }{\to }C\otimes C\stackrel {f\otimes g}{\to}M\otimes M\stackrel {\mu }{\to }M$.

\begin{proposition}
Для комоноида $C\in {K}_0$ и моноида $M\in {\cal K}_0$ $Hom$-множество ${\cal K}(C,M)$ является моноидом относительно операции свертки $*:{\cal K}(C,M)\times {\cal K}(C,M)\to {\cal K}(C,M):(f,g)\mapsto \mu \circ (f\otimes g)\circ \delta $, то есть операция $*$
\begin{itemize}
\item ассоциативна: $\forall f,g,h\in {\cal K}(C,M)$ $(f*g)*h=f*(g*h)$, и
\item имеет единицу, равную $\eta \circ \varepsilon :C\to M$: $\forall f\in {\cal K}(C,M)$ $(\eta \circ \varepsilon )*f=f=f*(\eta \circ \varepsilon )$.
\end{itemize}
\end{proposition} 
Доказательство.
\begin{itemize}
\item Рассмотрим коммутативную диаграмму \\
\vskip -0.35cm \hskip 1cm$\vcenter{\xymatrix{
C \ar[d]_-{\delta } \ar[rr]^-{\delta } & & C\otimes C \ar[d]^-{\delta \otimes 1_C} &  \\
C\otimes C \ar[r]_-{1_C\otimes \delta } & C\otimes (C\otimes C) \ar[r]^-{\alpha ^{-1}_{C,C,C}} \ar[d]_-{f\otimes (g\otimes h)} & (C\otimes C)\otimes C \ar[d]^-{(f\otimes g)\otimes h} &  \\
 &  A\times (A\otimes A) \ar[r]^-{\alpha ^{-1}_{A,A,A}} \ar[d]_-{1_A\otimes \mu } & (A\otimes A)\otimes A \ar[r]^-{\mu \otimes 1_A} & A\otimes A \ar[d]^-{\mu } \\
 & A\otimes A \ar[rr]_-{\mu } & & A
}}$

\noindent Верхний параллелограмм это условие того, что $\delta $ -- коумножение в $C$, соответственно нижний -- условие того, что $\mu $ -- умножение в $M$, а средний -- это условие естественности $\alpha $. Если написать равенство произведений стрелок, взятых вдоль верхнего и нижнего 'полупериметров', то получим ассоциативность $(f*g)*h=f*(g*h)$.  

\item Рассмотрим коммутативную диаграмму\\
\vskip -0.35cm \hskip 1cm$\vcenter{\xymatrix{
 & C \ar[dl]_-{\lambda ^{-1}_C} \ar[d]^-{\delta } && C \ar[dr]^-{\rho ^{-1}_C} \ar[d]_-{\delta } & \\
I\otimes C \ar[d]_-{1_I\otimes f} & C\otimes C \ar[l]^-{\varepsilon \otimes 1_C} \ar[d]^-{(\eta \circ \varepsilon )\otimes f} && C\otimes C \ar[r]_-{1_C\otimes \varepsilon } \ar[d]_-{f\otimes (\eta \circ \varepsilon )} & C\otimes I \ar[d]^-{f\otimes 1_I} \\
I\otimes M \ar[r]^-{\eta \otimes 1_M} \ar[dr]_-{\lambda _M} & M\otimes M \ar[d]^-{\mu } && M\otimes M \ar[d]_-{\mu } & M\otimes I \ar[l]_-{1_M\otimes \eta } \ar[dl]^-{\rho _M} \\
 & M && M &  
}}$

\noindent Верхние треугольники выражают условие того, что $\varepsilon $ -- коединица в $C$, нижние -- что $\eta $ -- единица в $M$, а средние прямоугольники коммутируют по построению. Тогда композиции стрелок, взятых по внешним путям от $C$ к $M$ выражают факт, что $\eta \circ \varepsilon $ -- левая и правая единица относительно умножения $*$, то есть $(\eta \circ \varepsilon )*f=f=f*(\eta \circ \varepsilon )$ (Замечание: $\lambda _M\circ (1_I\otimes f)\circ \lambda ^{-1}_C=f$ и $\rho _M\circ (f\otimes 1_I)\circ \rho ^{-1}_C=f$ в силу естественности $\lambda $ и $\rho $).    \hfill  $\square $
\end{itemize} 

\addvspace{0.1cm}
В силу предыдущей теоремы $Hom$-множество ${\cal K}(B,B)$ бимоноида $B$ является моноидом относительно свертки. 

\begin{definition}
Пусть $B\in {\cal K}_0$ бимоноид в ${\cal K}$. 
\begin{itemize}
\item {\bf Антиподом} на $B$ называется стрелка $S:B\to B$, являющаяся обратным элементом единичной стрелки $1_B:B\to B$ относительно операции свертки, то есть $S*1_B=\eta \circ \varepsilon =1_B*S$. 
\item Бимоноид $B$, снабженный антиподом $S$, называется {\bf алгеброй Хопфа}.
\end{itemize} 
\end{definition}

Не каждый бимоноид имеет антипод, но если он существует, то единственный.

\addvspace{0.2cm}
{\bf Примеры.}
\begin{enumerate}
\item Рассмотрим декартову тензорную категорию $(\mathbf {Set},\times ,\mathbf 1)$. Мы знаем, что любое множество $A$ имеет единственную структуру комоноида с коумножением $\delta :A\to A\times A:a\mapsto (a,a)$ и коединицей $!:A\to \mathbf 1$. Возьмем теперь множество $G$, являющееся группой (то есть моноидом, в котором каждый элемент имеет обратный). Структура коалгебры на $G$ совместима со структурой моноида на $G$ (с групповым умножением $\mu :G\times G\to G:(g_1,g_2)\mapsto g_1\cdot g_2$ и групповой единицей $\eta :\mathbf 1\to G:*\mapsto e$, где $*$ единственный элемент $\mathbf 1$, и $e$ единица $G$), то есть $G$ является бимоноидом в $(\mathbf {Set},\times ,\mathbf 1)$. Свертка в множестве функций $\mathbf {Set}(G,G)$ определяется как $(f*g)(h):=\mu \circ (f\times g)\circ \delta (h)=f(h)\cdot g(h)$, где $h\in G$. Соответственно антипод $S$ удовлетворяет равенству $e=e(h)=S(h)\cdot h=h\cdot S(h)$, то есть является операцией взятия обратного элемента и, поэтому, существует на группе. Если мы наоборот стартуем с определения алгебры Хопфа в $\mathbf {Set}$, то получим, что это должен быть моноид, (автоматически) совместимый с единственно возможной диагональной структурой комоноида, в котором каждый элемент имеет обратный по отношению к умножению, то есть придем к понятию группы. Таким образом, алгебры Хопфа в категории множеств $\mathbf {Set}$ совпадают с группами. 
\item Пусть $G$ -- группа с умножением $\cdot $ и единицей $e$. Рассмотрим векторное пространство $k(G)$ с базисом $G$ и коэффициентами из поля $k$. Тогда в тензорной категории векторных пространств $(\mathbf {Vect},\otimes ,k)$ пространство $k(G)$ является алгеброй Хопфа с умножением $\mu :k(G)\otimes k(G)\to k(G):g_1\otimes g_2\mapsto g_1\cdot g_2$, единицей $\eta :k\to k(G):1\mapsto e$, коумножением $\delta :k(G)\to k(G)\otimes k(G):g\mapsto g\otimes g$, коединицей $\varepsilon :k(G)\to k:g\mapsto 1$ и антиподом $S:k(G)\to k(G):g\mapsto g^{-1}$ (все отображения заданы на базисных элементах и продолжаются на другие по линейности).
\item Пусть $G$ -- конечная группа, состоящая из $n$ элементов. Рассмотрим множество функций на $G$ со значениями в поле $k$, то есть $\mathbf {Set}(G,k)$. Множество $\mathbf {Set}(G,k)$ имеет структуру векторного пространства размерности $n$ (действительно, возьмем в качестве базиса функции $\{f_g\, |\, g\in G\}$ такие, что $f_g(h):=\left\{
\begin{array}{lr}
1 & \text{если }h=g\\
0 & \text{если }h\ne g
\end{array}
\right.
$. Тогда любая функция $f:G\to k$ является линейной комбинацией функций $f_g$, именно, $f=\sum \limits _{g\in G}f(g)f_g$). Аналогично, $\mathbf {Set}(G\times G,k)$ является векторным пространством размерности $n^2$ с базисом $\{f_{(g,h)}\}$, $(g,h)\in G\times G$. Отождествим базисные элементы пространств $\mathbf {Set}(G\times G,k)$ и $\mathbf {Set}(G,k)\otimes \mathbf {Set}(G,k)$ по правилу $f_{(g,h)}\leftrightarrow f_g\otimes f_h$, то есть получим канонический изоморфизм этих пространств. Векторное пространство $\mathbf {Set}(G,k)$ имеет структуру алгебры Хопфа в категории $(\mathbf {Vect},\otimes ,k)$. Операции задаются следующим образом: умножение $\mu :\mathbf {Set}(G,k)\otimes \mathbf {Set}(G,k)\to \mathbf {Set}(G,k):f_{g_1}\otimes f_{g_2}\mapsto \left\{
\begin{array}{lr}
f_{g_1} & \text{если }g_1=g_2\\ 
0 & \text{если }g_1\ne g_2
\end{array}
\right. $,
единица $\eta :k\to \mathbf {Set}(G,k):1\mapsto \sum \limits _{g\in G}1\cdot f_g$,
коумножение $\varepsilon :\mathbf {Set}(G,k)\to \mathbf {Set}(G,k)\otimes \mathbf {Set}(G,k):(f=\sum \limits _{g\in G}f(g)f_g)\mapsto (f\circ \mu =\sum \limits _{g\in G}f(g)(\sum \limits _{\stackrel {\scriptstyle g_1,g_2\in G}{g_1\cdot g_2=g}}f_{g_1}\otimes f_{g_2}))$,
коединица $\varepsilon :\mathbf {Set}(G,k)\to k:f\mapsto f(e)$ (где $e$ -- единица группы $G$),
антипод $S:\mathbf {Set}\to \mathbf {Set}:f\mapsto f\circ i$ (где $i:G\to G:g\mapsto g^{-1}$ -- операция обращения в группе).
\end{enumerate}

Пример 3 является наиболее частым прототипом 'реальных' алгебр Хопфа, естественно возникающих в геометрии, комбинаторике, квантовой физике, хотя наверное еще раньше они возникли в топологии (в теории гомотопий) в абстрактном виде, близком к нашему определению. 

\addvspace{0.2cm}
В заключение приведем одно утверждение без доказательства.

\begin{proposition}
\begin{itemize}
\item Антипод коммутирует с гомоморфизмами бимоноидов, то есть если $f:B_1\to B_2$ такой гомоморфизм, и $S_i:B_i\to B_i$, $i=1,2$ -- антиподы в $B_1$ и $B_2$, то $f\circ S_1=S_2\circ f$.
\item В силу предыдущего пункта категория алгебр Хопфа является {\bf полной} подкатегорией категории бимоноидов.
\end{itemize}
\end{proposition}

\begin{assignment}
\begin{enumerate}
\item Убедиться, что в определении бимоноида (определение 21) оба условия ($\mu $ и $\eta $ -- гомоморфизмы комоноидов, и $\delta $ и $\varepsilon $ -- гомоморфизмы моноидов) эквивалентны.
\item Разобраться в деталях доказательства утверждения 16.
\item Доказать, что антипод, если существует, то единственный.
\item Доказать, что структура комоноида совместима со структурой моноида на группе $G$ в $(\mathbf {Set},\times ,\mathbf 1)$ (пример 1).
\item Проверить, что векторное пространство $k(G)$ (пример 2) действительно является алгеброй Хопфа.
\item Разобраться в структуре операций алгебры Хопфа $\mathbf {Set}(G,k)$ (пример 3) и проверить выполнение аксиом.
\end{enumerate}
\end{assignment}

\end{document}